\documentclass[11pt]{amsart}
\usepackage{amsmath,amssymb,amscd}
\usepackage[T1]{fontenc}
\usepackage{lmodern}
\usepackage{url}
\usepackage[hidelinks,bookmarksnumbered,unicode]{hyperref}

\newtheorem{theorem}{Theorem}[section]
\newtheorem{proposition}[theorem]{Proposition}
\newtheorem{lemma}[theorem]{Lemma}
\newtheorem{corollary}[theorem]{Corollary}
\theoremstyle{definition}
\newtheorem{quest}[theorem]{Question}
\theoremstyle{remark}

\DeclareMathOperator{\Pic}{Pic}
\DeclareMathOperator{\Hom}{Hom}
\newcommand{\OO}{\mathcal O}
\newcommand{\PP}{\mathbb P}
\newcommand{\QQ}{\mathbb Q}
\newcommand{\RR}{\mathbb R}
\newcommand{\XX}{\mathcal X}
\newcommand{\TT}{\mathcal T}

\newcommand{\DD}{\Delta}

\newcommand{\simq}{\sim_{\QQ}}

\title[Smooth limits with vanishing first Chern class]{Smooth limits of compact K\"ahler manifolds\\ with vanishing first Chern class, II}
\author{Mu-Lin Li}
\address{School of Mathematics, Hunan University, China}
\email{mulin@hnu.edu.cn}
\date{}
\begin{document}
\begin{abstract}
Let \(X\) be a compact complex manifold of dimension \(n\) that is deformation equivalent to a compact K\"ahler manifold with vanishing real first Chern class. We prove that \(X\) is K\"ahler and has vanishing real first Chern class when \(n\le3\). For every \(n\ge4\) we construct a smooth family of compact complex manifolds over the disc whose fibers over the punctured disc are compact K\"ahler manifolds with vanishing real first Chern class, while the central fiber is a non-K\"ahler Moishezon manifold.
\end{abstract}
\maketitle

\section{Introduction}

Let \(\DD\) be a complex disc centered at \(0\), and set
\(\DD^*=\DD\setminus\{0\}\). Two compact complex manifolds \(X\) and
\(X'\) are \emph{deformation equivalent} if there exist a proper
holomorphic submersion \(f:\XX\to B\) over a connected complex manifold
\(B\) and points \(b,b'\in B\) with \(X_b\cong X\) and \(X_{b'}\cong X'\).
Equivalently, one may use a finite chain of disc families: there exist
compact complex manifolds \(X_0=X,X_1,\dots,X_k=X'\) and proper
holomorphic submersions \(\pi_i:\XX_i\to\DD\) having fibers isomorphic to
\(X_{i-1}\) and \(X_i\).

Yau asked whether the K\"ahler property persists under deformation for
K\"ahler--Einstein manifolds of nonpositive scalar curvature.

\begin{quest}\cite[Chapter~IX, \S~IV, Problem~55]{SY94}\label{quest:yau}
If a compact complex manifold \(X\) can be deformed to a
K\"ahler--Einstein manifold \(X'\) of nonpositive scalar curvature, is
\(X\) K\"ahler?
\end{quest}

Hironaka's example \cite{Hironaka} shows that a smooth limit of projective
manifolds need not be K\"ahler. Yau's question asks whether the Einstein
hypothesis rules out such behavior.

Combining Siu's theorem \cite{SiuK3} with the classification of surfaces,
one obtains that \(X\) is K\"ahler whenever it can be deformed to a smooth
K\"ahler surface. The author \cite{LiStable} answers the question
affirmatively in the case where \(X\) can be deformed to a
K\"ahler--Einstein manifold \(X'\) of negative scalar curvature.
The author \cite[Theorem~1.2]{LiSmooth} proves that a smooth limit of
compact K\"ahler manifolds with vanishing real first Chern class belongs
to Fujiki's class \(\mathcal C\), that is, is bimeromorphic to a compact
K\"ahler manifold. The present paper sharpens this conclusion in relative
dimension three.

\begin{theorem}[cf. Theorem \ref{thm:global-1}]\label{thm:global}
Let \(\pi:\XX\to\DD\) be a proper holomorphic submersion with connected
fibers of dimension \(n=3\). Suppose that for some \(t_*\in\DD\) the fiber
\(X_{t_*}\) is K\"ahler and
\[
c_1(X_{t_*})=0\quad\text{in }H^2(X_{t_*},\RR).
\]
Then for every \(t\in\DD\) the fiber \(X_t\) is a compact K\"ahler manifold
and \(c_1(X_t)=0\) in \(H^2(X_t,\RR)\).
\end{theorem}

Theorem~\ref{thm:global} thus answers the Ricci-flat case of Question~\ref{quest:yau} affirmatively in dimensions at most three. The dimension
bound is sharp, as the following consequence of the constructions of
Oguiso--Schr\"oer \cite{OguisoSchroer} and Huybrechts
\cite{HuybrechtsFlop} shows.

\begin{theorem}[Theorem \ref{prop:examples-even}]\label{thm:counterexamples}
For each integer \(n\ge4\) there is a proper holomorphic submersion with
connected \(n\)-dimensional fibers over \(\DD\) whose punctured fibers are
compact K\"ahler manifolds with trivial canonical bundle and whose central
fiber is non-K\"ahler. The central fiber is smooth and Moishezon with
trivial canonical bundle. In dimension four the punctured fibers can be
taken to be irreducible holomorphic symplectic.
\end{theorem}

The Beauville--Bogomolov decomposition reduces the threefold problem,
after a finite \'etale cover of one fiber, to complex tori, products of
an elliptic curve with a K3 surface, and simply connected Calabi--Yau
threefolds. In this paper a \emph{projective Calabi--Yau threefold}
is smooth and connected, with $K_X\simeq\OO_X$ and
$H^1(X,\OO_X)=H^2(X,\OO_X)=0$.  The main projectivity result is the following.

\begin{theorem}[cf. Theorem \ref{thm:projective-1}]\label{thm:projective}
Let $\pi:\XX\to\DD$ be a smooth family whose punctured fibers are
simply connected projective Calabi--Yau threefolds. Then $X_0$ is a
simply connected projective Calabi--Yau threefold.
\end{theorem}

Rao's specialization theorem first gives that $X_0$ is Moishezon.
A projective small model $T$ exists, but its ample classes can change
degree across the small modification. We compare $X_0$ and $T$
inside smooth families in Proposition~\ref{prop:comparison}.
For a punctured polarization $L$, write $L_0=P_0+N$ for its divisorial
Zariski decomposition. Lutz's description of the movable cone
identifies the possible fixed prime divisors. The surface index
theorem gives $(P_0^2N)_{X_0}\le0$, whereas a nonempty sequence
of $P_t$-negative flops would give the opposite strict inequality
by Wilson's cubic formula. This proves nefness of the marked class
$P_t$ on the given punctured fibers and yields a semiample $P_0$.
Its contraction $f:X_0\to Z$ has fibers of dimension at most one
(Theorem~\ref{thm:main}). When $b_2=1$, a direct comparison of two
polarizations already proves projectivity
(Corollary~\ref{prop:rankone-direct}).

For arbitrary $b_2$, a single integral class avoids the roots of
the Du Val slices over the general points of the exceptional image
curves. We run a local minimal model program directed by this class
over each point of $Z$. Each step preserves the central model over
the target (Lemma~\ref{lem:localflop}). A curve of degree zero on a
final model would deform to a curve in a general slice and violate
the choice of class (Lemma~\ref{lem:root}). The local ample models
are thus relative Proj spaces of one graded algebra, whose
tautological bundle polarizes $f$. The elliptic-curve times K3 case admits a different argument. Periods of a relative one-form and a holomorphic vector field show that the product structure persists under smooth deformation (Lemma~\ref{lem:marked-limits} and Proposition~\ref{prop:genericK}). Section~\ref{sec:propagation} gives the one-fiber theorem by isolating the exceptional parameters on the disc.

\section{Preliminaries}\label{sec:preliminaries}

\subsection{Smooth families and finite covers}\label{subsec:bb}
Write $X_t=\pi^{-1}(t)$ and $K_{\XX/\DD}=\det\Omega^1_{\XX/\DD}$.
Ehresmann's theorem identifies the differentiable fibers over a disc.
We use the resulting local systems to mark cohomology. A modification
is a proper bimeromorphic morphism. A compact complex manifold in
Fujiki's class $\mathcal C$ admits a modification from a compact
K\"ahler manifold and satisfies the $\partial\bar\partial$-lemma
\cite[Corollary~5.23]{DGMS}. In particular its holomorphic forms are
closed and its cohomology has a Hodge decomposition compatible with
conjugation and cup products. A Moishezon manifold is bimeromorphic
to a projective variety and has a smooth projective modification.

If all the fibers in a smooth family satisfy the
$\partial\bar\partial$-lemma, their Hodge numbers are locally
constant. Indeed, each $h^{p,q}$ is upper semicontinuous and their
sum for fixed $p+q=k$ equals the constant Betti number $b_k$.
A compact K\"ahler manifold with $h^{2,0}=0$ is projective, since its
open K\"ahler cone contains a rational class in $H^2(X,\RR)$.

A compact K\"ahler manifold with $c_1(X)=0$ in $H^2(X,\RR)$ admits a
connected finite \'etale cover that is a product of a torus and simply
connected factors with trivial canonical bundle
\cite[Theorem~2]{BeauvilleKummer}. A strict Calabi--Yau factor has no
nonzero holomorphic forms in intermediate degrees. An irreducible
holomorphic symplectic factor has even dimension and its holomorphic
two-forms are spanned by a nondegenerate form. In dimension three the
possibilities for the cover are a torus, a product of an elliptic
curve and a K3 surface, and a simply connected strict Calabi--Yau
threefold. In the last case $h^{2,0}=0$ gives projectivity. A K3
surface is K\"ahler \cite{SiuK3}. K\"ahlerness descends under a finite
\'etale map by taking the trace of a K\"ahler form in local inverse
branches.

\subsection{Holomorphic forms and ampleness in families}\label{subsec:propinputs}
K\"ahlerness is open under small deformations \cite{KodairaSpencer}.
On a K\"ahler neighborhood of a parameter $s$ the Hodge numbers are
constant. Grauert's theorem identifies the fiber of
$\pi_*\Omega^p_{\XX/\DD}$ at $s$ with $H^0(X_s,\Omega^p_{X_s})$
\cite{Grauert}. Cartan's theorem~A on the disc then extends any
prescribed holomorphic $p$-form on $X_s$ to a global relative form.
The extension is relatively closed by the identity theorem. Rao proved that a smooth limit of Moishezon manifolds is Moishezon.
Rao and the coauthors of the appendix also showed that the projective
locus in a smooth disc family with uncountably many projective fibers
is the complement of a closed discrete set
\cite[Corollaries~1.4 and~A.5]{Rao}.

\begin{lemma}\cite[Chapter~IV, Theorem~4.2]{BS76}\label{lem:amplelocus}
Let $L$ be a line bundle on a smooth family $\pi:\XX\to\DD$.
If $L_s$ is ample for one parameter $s$, then $L_t$ is ample
outside a closed discrete subset of $\DD$.
\end{lemma}

\subsection{Divisors and small modifications}\label{subsec:positivity}
We use additive notation in $\Pic(X)\otimes\QQ$. On an $n$-fold a
rational divisor $D$ is big if its volume is positive, where the
volume is the limsup of $n!h^0(X,mD)/m^n$ for sufficiently divisible
$m$. We call $D$ movable if some multiple has a nonempty complete
linear system without a fixed prime divisor. On a projective model
$\overline{\operatorname{Mov}}$ denotes the closure of the cone these
classes generate. A semiample class has a globally generated multiple.
For a big semiample class $P$ the associated Stein morphism $g$ is
birational, its target is normal projective, and $P=g^*H$ with $H$
ample and rational. For a finitely generated big section ring the
fixed divisor of a Veronese generated in degree one grows linearly.
Subtracting it gives the divisorial Zariski decomposition used below.

A small modification is a bimeromorphic map inducing an isomorphism
on complements of analytic sets of codimension at least two on
both sides. A proper bimeromorphic morphism is small if it contracts
no divisor. On normal spaces, strict transform across a small
modification identifies prime divisors and meromorphic sections with
a given divisor bound. In particular it identifies complete section
rings. For models over a common base it also identifies their relative
section algebras. If a line bundle is relatively ample, the relative
Proj of its section algebra recovers its source. For a proper
analytic morphism, fiberwise ampleness implies relative ampleness
locally on the target \cite[Remark~3.2]{FujinoPositive}.

\subsection{Minimal models and the flop formula}\label{subsec:birational}
We distinguish algebraic $\QQ$-factoriality on projective varieties
from $\QQ$-factoriality for local analytic Weil divisors. On a
resolution $p:W\to Z$, write
$K_W=p^*K_Z+\sum a(E,Z)E$. The singularities are canonical if the
exceptional discrepancies are nonnegative, and terminal if they are
positive. Canonical singularities are rational. A crepant morphism
preserves the canonical divisor under pullback. We use the discrepancy
conventions for klt and plt pairs in \cite[Chapter~5]{KMbook}.
In particular the coefficient-one boundary of a plt pair is normal.

In the projective threefold arguments a big effective directing class
$D$ on a smooth $K$-trivial model is treated by running the minimal
model program (MMP) for $(X,\varepsilon D)$ with $\varepsilon>0$
sufficiently small. At a step with scaling parameter $\lambda>0$,
the scaled adjoint class vanishes on the selected ray. For an elementary
flop, choose an effective generator $\gamma$ on its source. On smooth
Calabi--Yau threefolds Wilson's formula gives, for corresponding
divisors $D,D'$,
\[
 (D')^3-D^3=-c(D\cdot\gamma)^3,\qquad c>0
 \quad\text{\cite[Proposition~2.1 and (2.1.1)]{WilsonFlops}}.
\]
The coefficient depends on the contraction and the normalization of
$\gamma$. The local MMP in Section~\ref{sec:local} is instead run for
a projective analytic morphism over a Stein neighborhood of a point.
Numerical equivalence there tests degrees against line bundles defined
near the whole fiber \cite[\S4]{FujinoMMP}.

\subsection{Conic fibrations and roots}\label{subsec:roots}
In the terminology of Lutz, a conic fibration surface is an irreducible
projective Cohen--Macaulay surface with a flat morphism with connected
fibers to a smooth curve, whose general fiber is a projective line or
two projective lines meeting in an ordinary node. An SQM-conic fibration
surface is a divisor that admits such a fibration on a small
$\QQ$-factorial modification. If the base has positive genus, the class
$\ell_E$ of the whole general conic is intrinsic to the marked divisor
$E$. It satisfies $E\cdot\ell_E=-2$. The reflection
$s_E(D)=D+(D\cdot\ell_E)E$ fixes $\ell_E^\perp$ and sends $E$ to
$-E$ \cite[Definition~1.2 and Lemma~3.8]{Lutz}.
We apply Lutz's cone comparison on the Kuranishi family of a
\emph{smooth projective} small model. Proposition~\ref{prop:comparison}
relates its marking to that of the original disc family.

A normal canonical surface singularity is Du Val. The irreducible
exceptional curves of its minimal resolution give a basis of the
second homology of a resolved small ball. Their intersection form is
the negative Cartan form of type $A$, $D$ or $E$
\cite[\S2.1]{PlSt}. A root is an integral vector of square $-2$.
There are finitely many roots and each has coefficients of one sign
in the curve basis. Indeed, in a decomposition $v=v_+-v_-$ into
nonnegative parts of disjoint support, both nonzero, one has
$v_+\cdot v_-\ge0$ and $v_\pm^2\le-2$ by evenness and negative
definiteness. Thus $v^2\le-4$. In Section~\ref{sec:closure} we also
verify that the images of these roots in global homology are nonzero.

\section{The family and its projective small models}\label{sec:family}

In this section we fix a proper holomorphic submersion
\(\pi:\XX\to\DD\) with connected three-dimensional fibers such that
\(X_t\) is a simply connected projective Calabi--Yau threefold for every
\(t\ne0\). We shrink the disc about zero whenever necessary, and we use
the marking supplied by Ehresmann's theorem.

\begin{lemma}\label{lem:family}
Let \(\pi:\XX\to\DD\) be a proper holomorphic submersion with connected three-dimensional fibers. Suppose that \(X_t\) is a simply connected projective Calabi--Yau threefold for every \(t\ne0\). Then, after shrinking \(\DD\) about the origin, the following hold.
\begin{enumerate}
\item The central fiber is Moishezon and simply connected. The relative canonical bundle is trivial, and \(H^1(X_t,\OO)=H^2(X_t,\OO)=0\) for all \(t\).
\item Under the topological marking, the restriction maps identify \(\Pic(\XX)\) and \(\Pic(X_t)\) with \(H^2(X_t,\mathbb Z)\). There is a line bundle \(L\) on \(\XX\) such that \(L_t\) is ample for all sufficiently small \(t\ne0\) and \(L_0\) is big.
\end{enumerate}
\end{lemma}

\begin{proof}
By \cite[Corollary~1.4]{Rao}, \(X_0\) is Moishezon, and by Ehresmann's
theorem it is simply connected. The coherent sheaf
\(\pi_*K_{\XX/\DD}\) is torsion-free of rank one, hence a line
bundle on the disc. A local generator near zero is a relative volume form that does not
vanish over the puncture, so its divisor is \(aX_0\) for some \(a\ge0\):
the central fiber is smooth, connected and irreducible. Dividing by
\(t^a\) therefore gives a nowhere vanishing relative canonical section,
and \(K_{\XX/\DD}\simeq\OO_{\XX}\) after shrinking.

The \(\partial\bar\partial\)-lemma on \(X_0\), together with
\(b_1(X_0)=0\), gives \(H^1(X_0,\OO)=0\). Serre duality and
\(K_{X_0}\simeq\OO_{X_0}\) give \(H^2(X_0,\OO)=0\).
The corresponding vanishings hold on the punctured fibers by
hypothesis. Grauert's theorem therefore gives \(R^j\pi_*\OO_{\XX}=0\) for \(j=1,2\),
and the Stein property of the disc together with the Leray spectral
sequence implies \(H^j(\XX,\OO_{\XX})=0\) for \(j=1,2\). Thus the exponential
sequence and Ehresmann's theorem identify
\[
 \Pic(\XX)\xrightarrow[\simeq]{c_1}H^2(\XX,\mathbb Z)
 \xrightarrow{\simeq}H^2(X_t,\mathbb Z)
 \xleftarrow[\simeq]{c_1}\Pic(X_t).
\]
These identifications are compatible with restriction.

Extend an ample bundle on one punctured fiber to a line bundle \(L\) on
\(\XX\). Lemma~\ref{lem:amplelocus} then shows that, after shrinking
about zero, \(L_t\) is ample for every \(t\ne0\).
For each \(k>0\), Kodaira vanishing gives
\(h^0(X_t,L_t^k)=\chi(X_t,L_t^k)\) on this punctured disc, and upper
semicontinuity together with constancy of the Euler characteristic gives
\[
 h^0(X_0,L_0^k)\ge\chi(X_t,L_t^k)
 =\frac{L_t^3}{6}k^3+O(k^2),
\]
where \(L_t^3>0\). Hence \(L_0\) is big.
\end{proof}

\begin{lemma}[Movement of a rational curve]\label{lem:movement}
Let \(\pi:\XX\to\DD\) be a smooth proper family of threefolds with
\(K_{X_0}\simeq\OO_{X_0}\). Let \(C\subset X_0\) be an integral rational
curve and \(M\in\Pic(\XX)\otimes\QQ\). Suppose either that \(M_t\) is ample
on a dense subset of \(\DD^*\) and \(M_0\cdot C\le0\), or that \(M_t\) is
nef on a dense subset of \(\DD^*\) and \(M_0\cdot C<0\). The normalization
map of \(C\) then deforms in \(X_0\) with infinitely many distinct image
curves, all of the same integral homology class as \(C\).
\end{lemma}

\begin{proof}
Write \(u:\PP^1\to\XX\) for the normalization map. The analytic space of
maps from the fixed source \(\PP^1\) has tangent space
\(H^0(\PP^1,u^*T_{\XX})\) and obstruction space
\(H^1(\PP^1,u^*T_{\XX})\). More explicitly, the linearization of the
holomorphic-map equation at \(u\) is the Dolbeault operator on
\(u^*T_{\XX}\). Splitting off its image and applying the implicit
function theorem gives the local Kuranishi equations on its kernel,
with values in its cokernel. Thus the reduced map germ has an
irreducible component \(B\) through \(u\) of complex dimension at
least \(\chi(\PP^1,u^*T_{\XX})\) (see also
\cite[Lemma~5.1]{LLR26}). The exact sequence
\[
 0\longrightarrow T_{X_0}\longrightarrow T_{\XX}|_{X_0}
 \longrightarrow \OO_{X_0}\longrightarrow0
\]
gives \(\operatorname{rank}(u^*T_{\XX})=4\) and
\(\deg(u^*T_{\XX})=0\). Hence \(\dim B\ge4\).

After shrinking \(B\), all its maps are birational onto their images.
Choose \(p\in\PP^1\) over a smooth point of \(C\) and a coordinate
disc \(D\ni p\) on which a target coordinate composed with \(u\) is
univalent on a neighborhood of \(\overline D\). For \(v\) close to
\(u\) that coordinate stays univalent on \(D\), by convergence in
\(C^1\) in a coordinate in which the original map is the identity.
Since \(u(p)\notin u(\PP^1\setminus D)\), compactness also gives
\(v(p)\notin v(\PP^1\setminus D)\). Thus \(v(p)\) has one unramified
preimage. Maps with a given image now differ only by
an automorphism of \(\PP^1\). Their germs have dimension at most three:
three marked points sent to distinct smooth points of the fixed image
recover the automorphism in local coordinates. Thus a component of
dimension at least four contains infinitely many different images.

For each \(v\in B\), the composite \(\pi\circ v\) is constant, and
evaluation at any fixed source point defines a holomorphic map
\(b:B\to\DD\). The pushforward of \([\PP^1]\) is locally constant in
\(H_2(\XX,\mathbb Z)\), as is its pairing with \(c_1(M)\).
If \(b\) were nonconstant, restrict to a normalized analytic curve
germ through \(u\) on which it is nonconstant. The open mapping
theorem makes its image meet the given dense set. The degree there
would be positive in the ample case and nonnegative in the nef case,
contrary to its value on \(C\). Thus \(b=0\) on the germ, and the
marking identifies all the resulting curve classes with \([C]\).
\end{proof}

\begin{proposition}\label{prop:comparison}
Let \(\pi:\XX\to\DD\) be a proper holomorphic submersion with connected
three-dimensional fibers. Assume that \(X_t\) is a simply connected
projective Calabi--Yau threefold for every \(t\ne0\).
\begin{enumerate}
\item There exist an algebraically \(\QQ\)-factorial normal projective
threefold \(T\) with terminal singularities and a bimeromorphic map
\(\phi:X_0\dashrightarrow T\) that is an isomorphism in codimension one.
\item For every pair \((T,\phi)\) as in~(1), \(T\) is smooth and
\(K_T\simeq\OO_T\). After shrinking \(\DD\) about zero there exist a
smooth projective morphism \(\rho:\TT\to\DD\), an identification
\(T_0\simeq T\), and a bimeromorphic map \(\Phi:\XX\dashrightarrow\TT\)
over \(\DD\). There are closed analytic subsets \(E_X\subset\XX\)
and \(E_T\subset\TT\) whose intersections with every fiber have
complex codimension at least two such that \(\Phi\) induces an
isomorphism \(\XX\setminus E_X\simeq\TT\setminus E_T\) over \(\DD\).
The induced maps \(\phi_t:X_t\dashrightarrow T_t\) are small, and
\(\phi_0=\phi\) under the given identification.
\item For the extension in~(2), restriction to the common open sets
defines an isomorphism of integral local systems
\[
\theta:R^2\pi_*\mathbb Z\xrightarrow{\ \simeq\ }R^2\rho_*\mathbb Z.
\]
For each \(t\), if \(M\) is a line bundle on \(X_t\) and \(M'\) is its
strict transform on \(T_t\), then \(\theta_t(c_1(M))=c_1(M')\). In
particular, the strict transforms of the fiber restrictions of a
line bundle on \(\XX\) are the restrictions of a single line bundle on
\(\TT\).
\end{enumerate}
\end{proposition}

\begin{proof}
\smallskip\noindent\emph{Step 1. Existence of a projective small model.}
By Lemma~\ref{lem:family}, the central fiber is simply connected and
Moishezon, with \(K_{X_0}\simeq\OO_{X_0}\) and
\(H^j(X_0,\OO_{X_0})=0\) for \(j=1,2\). Take a smooth projective
modification of \(X_0\). The projective threefold minimal model
program and abundance produce an algebraically \(\QQ\)-factorial
projective terminal good minimal model \(T\) with \(K_T\simq0\)
\cite[proof of Theorem~3.8, Step~1]{BCDG}.

Let \(\eta\) be a nowhere vanishing holomorphic three-form on
\(X_0\), and denote by \(\eta_T\) the induced rational canonical
form on \(T\). On a smooth common resolution
\[
 a:W\longrightarrow X_0,\qquad b:W\longrightarrow T,
\]
the form \(a^*\eta\) is holomorphic. Thus \(\eta_T\) has no pole
along a prime divisor of \(T\), so it is a section of the reflexive
canonical sheaf. Its canonical Weil divisor is effective and
\(\QQ\)-linearly trivial, hence zero. Indeed, a positive multiple
would be the effective divisor of a rational function on the normal
projective variety \(T\), and such a function is constant. We may
therefore choose \(K_T=0\) as a Cartier divisor.

The divisor of the rational form \(a^*\eta=b^*\eta_T\) on \(W\) is
supported precisely on the \(a\)-exceptional prime divisors, because
\(X_0\) is smooth. Terminality of \(T\) identifies its support with
the union of the \(b\)-exceptional prime divisors. These two sets therefore
coincide, so the induced map
\(X_0\dashrightarrow T\) is an isomorphism in codimension one.
This proves~(1).

\smallskip\noindent\emph{Step 2. A smooth analytic
\(\QQ\)-factorialization of the prescribed model.}
Fix any pair \((T,\phi)\) satisfying~(1). Smallness transfers the form \(\eta\) to a rational canonical form with
zero Weil divisor on this \(T\) as well. Thus \(K_T=0\) is Cartier and
\(\omega_T\simeq\OO_T\).

The model \(T\) is at this stage only algebraically
\(\QQ\)-factorial. Take an analytic \(\QQ\)-factorialization
\cite[Lemma~7.8]{IS},
\[
 \tau_0:U\longrightarrow T.
\]
This is a proper small bimeromorphic morphism from an analytically
\(\QQ\)-factorial compact normal threefold. Smallness and \(K_T=0\) give
\(K_U=\tau_0^*K_T=0\), and \(U\) is terminal. The induced map
\[
 \psi_0:X_0\dashrightarrow U
\]
is now between compact \emph{analytically} \(\QQ\)-factorial
terminal threefolds with nef canonical divisors. The analytic form of
\cite[Theorem~4.9]{Kflops} factors it into finitely many analytic
flops. By \cite[Corollary~4.11]{Kflops}, these flops preserve the
analytic singularities. Since \(X_0\) is smooth, all intermediate
models, including \(U\), are smooth.

\smallskip\noindent\emph{Step 3. Extending the analytic flops.}
Consider an elementary flop in the factorization from Step~2, with
central contraction diagram
\[
 V^-_0\xrightarrow{f_0} W_0\xleftarrow{f^+_0} V^+_0 .
\]
Suppose we have already constructed the deformation
\(\mathcal V^-\to\DD\). All the central source models are smooth
with trivial canonical bundle. The contraction target is normal and
the transferred volume form has zero canonical Weil divisor.
Consequently \(W_0\) is Gorenstein and terminal by the discrepancy
comparison for a small map from a smooth threefold. In particular
\(Rf_{0*}\OO_{V^-_0}=\OO_{W_0}\), and the exceptional locus is compact. The analytic blowdown deformation
theorem \cite[Lemma~35]{KollarBlowdown}, cf.\ \cite[Proposition~11.4]{KMflips}, extends \(f_0\) to a
proper morphism
\[
 F^-:\mathcal V^-\longrightarrow\mathcal W
\]
over a smaller disc, with \(\mathcal W\) a flat deformation of
\(W_0\). After shrinking, \(\mathcal W\) and its fibers are normal,
and the fibers are irreducible. On an open set of the central
isomorphism locus the inverse function theorem and properness give a
one-sheeted open set for the nearby maps. Hence \(F^-\) and its
fiber maps are bimeromorphic there and therefore bimeromorphic
onto their targets.

We describe the opposite family explicitly. The exceptional images
of \(f_0\) are finitely many points. At any such point \(w\), the
terminal Gorenstein threefold germ is an isolated compound Du Val
hypersurface double point. There are coordinates in which its
equation is
\[
 x^2+g(y_1,y_2,y_3)=0 .
\]
In these coordinates Koll\'ar's local construction realizes
\(f_0^+\) by composing \(f_0\) with the deck involution \(x\mapsto-x\)
\cite[Section~2.2, Example~2.3 and Theorem~2.4]{Kflops}.
Indeed, for a Weil divisor \(D\) on the double cover, the sum of
\(D\) and its conjugate is the pullback of a Weil divisor on its
smooth quotient. Thus the involution changes the sign of the local
divisor class and gives the prescribed opposite contraction.

Lift the four coordinates to the flat deformation germ. Together
with the disc coordinate \(t\) they realize it as a hypersurface
in \(\mathbb C^4\times\DD\). Indeed, the lifted coordinates give a local
embedding by Nakayama, and the central hypersurface equation lifts
to a generator of the ideal by flatness and Nakayama. Weierstrass
preparation in \(x\) writes the equation, up to a unit, as
\[
 x^2+a(y,t)x+b(y,t)=0,\qquad a(y,0)=0,\quad b(y,0)=g(y).
\]
The involution \(\iota(x,y,t)=(-a(y,t)-x,y,t)\) thus extends the
specified central involution over the original disc.

Take disjoint invariant neighborhoods of the exceptional images in
\(\mathcal W\). On each, use a copy of the inverse image under \(F^-\)
with its map to that neighborhood given by \(\iota\circ F^-\).
Write \(B^-\) for the closed locus of positive-dimensional fibers of
\(F^-\). By properness, after shrinking the disc, \(B^-\) is contained
in relatively smaller neighborhoods of these same points. The union
\(B\) of \(B^-\) and its local images under the involutions is
therefore closed analytic. On \(\mathcal W\setminus B\), the original
maps and the modified local maps are isomorphisms, by normality and
the finite mapping theorem. Glue the modified local pieces to
\(\mathcal W\setminus B\) using these isomorphisms over
\(\mathcal W\). The resulting morphism
\(F^+:\mathcal V^+\to\mathcal W\) is proper, since properness is
local on the target. Its source is smooth over the disc on each
chart. Its central restriction is the prescribed flop with its given
identification on the common open set. This is the Gorenstein case
of the flop deformation construction of
\cite[Theorem~11.10]{KMflips}.

The inverse images \((F^\pm)^{-1}(B)\) have dimension at most one on
the central fibers. They are proper over the disc, so the same bound
holds on all nearby fibers. Off these sets the two families are
isomorphic over the disc. For a composite of two such identifications,
intersect their common open sets on the middle space. The complements
on the outer spaces are analytic. Each is the union of the old
complement and the image of the additional excluded set under the
proper graph correspondence. Outside the old complement
this is the inverse image under an isomorphism of a set of fiber
dimension at most one. Thus the bound persists under composition.
Induction gives a smooth proper family
\(p:\mathcal U\to\DD\) with \(U_0=U\) and a bimeromorphic map
\(\Psi:\XX\dashrightarrow\mathcal U\) extending \(\psi_0\), with the
claimed codimension bound on every fiber.

A common resolution of \(U\) and \(X_0\) gives
\(H^j(U,\OO_U)=0\) for \(j=1,2\). Upper semicontinuity gives the
same vanishings on every nearby \(U_t\). The exponential sequence,
Grauert's theorem, the Stein property of the disc and Ehresmann's
theorem therefore give restriction-compatible isomorphisms
\[
 \Pic(\mathcal U)\xrightarrow[\simeq]{c_1}
 H^2(\mathcal U,\mathbb Z)\xrightarrow{\simeq}H^2(U_t,\mathbb Z),
 \qquad
 \Pic(U_t)\xrightarrow[\simeq]{c_1}H^2(U_t,\mathbb Z).
\]

\smallskip\noindent\emph{Step 4. Constructing the projective family.}
The classes pulled back from \(T\) span \(H^2(U,\QQ)\).
Indeed every line bundle on the Moishezon manifold \(U\) has a meromorphic
section, by pullback to a smooth projective modification and descent
of a meromorphic section. For its divisor \(D\), the pushforward
\(D_T=(\tau_0)_*D\) is algebraic by Chow's theorem. Algebraic
\(\QQ\)-factoriality makes \(mD_T\) Cartier for some \(m>0\), and
smallness gives \(\tau_0^*(mD_T)=mD\). The assertion follows from the
Picard isomorphism of Step~3.

Choose very ample bundles \(A_0,\ldots,A_r\) on \(T\) whose pullbacks
span \(H^2(U,\QQ)\) and for which \(H^q(T,A_i)=0\) for \(q>0\).
One obtains such bundles by twisting a finite spanning set by high
powers of one ample bundle and adjoining a high power of that bundle.
Rational singularities and the projection formula give
\[
 H^q(U,\tau_0^*A_i)=H^q(T,A_i)=0 \quad(q>0).
\]
Extend each \(\tau_0^*A_i\) to a line bundle \(\mathcal L_i\) on
\(\mathcal U\). After shrinking, cohomology and base change and the
surjectivity of evaluation on \(U\) make these bundles relatively
generated. They define a proper map
\[
 h:\mathcal U\longrightarrow\mathcal P=
 \prod_{\DD,\,0\le i\le r}\mathbb P_{\DD}(p_*\mathcal L_i).
\]
Its central restriction is \(i\circ\tau_0\), where \(i:T\hookrightarrow
\mathcal P_0\) is a closed embedding. Thus
\((h_0)_*\OO_U=i_*\OO_T\) and \(R^1(h_0)_*\OO_U=0\).

Write \(\mathcal A=h_*\OO_{\mathcal U}\). The exact sequence for
multiplication by \(t\), followed by Nakayama's lemma and properness
of \(\mathcal P\) over the disc, gives \(R^1h_*\OO_{\mathcal U}=0\)
after shrinking. Using multiplication by \(t-a\) now gives
\[
 \mathcal A/(t-a)\mathcal A=(h_a)_*\OO_{U_a}\qquad(a\in\DD).
\]
The unit \(\OO_{\mathcal P}\to\mathcal A\) is surjective on the central
fiber, hence near that fiber. The Stein factorization is therefore
\[
 \mathcal U\xrightarrow{\tau}\TT\hookrightarrow\mathcal P .
\]
Here \(\TT\) is normal, and it is flat over the disc because
\(\mathcal A\) has no base torsion. It is projective over the disc,
with \(T_0=T\) scheme-theoretically. By the displayed base change
each \(T_a\) is the Stein target of \(h_a\), hence normal and
irreducible, and \((\tau_a)_*\OO_{U_a}=\OO_{T_a}\).
Let \(E\) be the inverse image under \(\tau\) of its locus of
positive-dimensional fibers. This is a closed analytic set proper
over \(\DD\). Since \(\dim(E\cap U)\le1\), the same bound holds
on all nearby fibers. Therefore \(\tau_a\) is finite away from
\(\tau_a(E\cap U_a)\), where it is an isomorphism by the structure
sheaf identity. In particular all the \(\tau_a\) are small
bimeromorphic morphisms.

\smallskip\noindent\emph{Step 5. The punctured fibers.}
Let \(A_{i,a}\) be the restriction to \(T_a\) of the tautological
bundle from the \(i\)-th projective factor. The equality
\(\mathcal L_i|_{U_a}=\tau_a^*A_{i,a}\) shows that the pullback
classes span \(H^2(U_a,\QQ)\), by the marking of \(\mathcal U\).
For an algebraic Weil divisor \(D\) on \(T_a\), its strict transform
\(\widetilde D\) on \(U_a\) is Cartier. The spanning property and the
Picard isomorphism of Step~3 give, after clearing denominators and
torsion,
\[
 \OO_{U_a}(m\widetilde D)
 \simeq\tau_a^*\!\left(\bigotimes_i A_{i,a}^{n_i}\right)
 \qquad (m>0,\ n_i\in\mathbb Z).
\]
The indicated bundle on \(T_a\) is algebraic by GAGA. Pushing down
the meromorphic equivalence gives \(mD\) linearly equivalent to a
Cartier divisor, since \(\tau_a\) is small and meromorphic functions
on the projective normal variety \(T_a\) are rational. This proves
algebraic \(\QQ\)-factoriality of \(T_a\).

Through \(X_a\dashrightarrow U_a\to T_a\) the volume form defines the
zero canonical Weil divisor on \(T_a\). Thus \(K_{T_a}\) is Cartier
and \(\tau_a\) is crepant. Each exceptional divisor over \(T_a\) is
also exceptional over the smooth \(U_a\), by smallness, and hence
has positive discrepancy. Thus \(T_a\) is terminal. For \(a\ne0\),
the comparison \(X_a\dashrightarrow T_a\) is between projective
algebraically \(\QQ\)-factorial terminal minimal threefolds.
The algebraic form of \cite[Theorem~4.9 and Corollary~4.11]{Kflops}
shows that \(T_a\) is smooth. Now \(\tau_a\) is a small morphism
between smooth threefolds. Its Jacobian determinant is nonzero off
a subset of codimension at least two, so it is nowhere zero. Thus
\(\tau_a\) is an isomorphism for \(a\ne0\).

\smallskip\noindent\emph{Step 6. The central fiber.}
A nontrivial \(\tau_0\) would have a positive-dimensional fiber.
Choose two distinct points of that fiber and lift them to a smooth
projective modification \(W\to U\). The composite \(W\to T\) is a
projective resolution of a klt variety, so a rational chain in its
fiber joins the chosen lifts \cite[Corollary~1.6]{HM}. Its image in
\(U\) contains a rational curve \(C\) contracted by \(\tau_0\).
For \(\mathcal L=\bigotimes_i\mathcal L_i\) one has
\[
  \mathcal L|_U=\tau_0^*\!\left(\bigotimes_i A_i\right),
  \qquad \mathcal L\cdot C=0 .
\]
Step~5 and the Segre polarization make \(\mathcal L|_{U_a}\) ample
for \(a\ne0\). Lemma~\ref{lem:movement} now yields infinitely many
distinct \(\mathcal L\)-null rational curves on \(U\). Every one is
contracted by \(\tau_0\), whereas its exceptional locus is a compact
analytic set of dimension at most one and thus contains only finitely
many irreducible curves. It follows that \(\tau_0\) is finite and hence
an isomorphism onto the normal \(T\). The morphism \(\tau\) is likewise
finite and bimeromorphic with normal target and so is an isomorphism.
Taking \(\Phi=\tau\circ\Psi\) proves~(2).

\smallskip\noindent\emph{Step 7. The marking.}
For a closed analytic subset \(A\) of codimension at least two in a
complex manifold \(V\), cohomology with supports \(H_A^i(V,\mathbb Z)\)
vanishes for \(i\le3\), by duality and the real dimension bound.
Thus removal of \(A\) preserves integral cohomology in degree two.
Apply this on the common open subsets of the two total spaces and
on each pair of fibers. The resulting identifications commute with
restriction. Since the disc is contractible, Ehresmann's theorem
identifies total degree-two cohomology with the sections of the
degree-two local system. This gives the asserted isomorphism
\(\theta\). Strict transforms agree on the common open sets, so
naturality of \(c_1\) identifies their Chern classes. Finally the
Picard isomorphisms on \(\mathcal U\simeq\TT\) and its fibers carry
each transferred integral class to the unique line bundle with these
restrictions. This proves~(3).
\end{proof}

\section{Specialization of a punctured polarization}\label{sec:polarization}

Throughout this section the punctured fibers are simply connected
projective Calabi--Yau threefolds. We use the Picard identifications
of Lemma~\ref{lem:family} and the comparison families of
Proposition~\ref{prop:comparison}.

\subsection{Two semiampleness criteria}

\begin{lemma}\label{lem:semiample}
Let \(M\) be a big movable rational line-bundle class on \(X_0\). If \(M\) has nonnegative degree on every rational curve in \(X_0\), then \(M\) is semiample, and its birational contraction is crepant.
\end{lemma}

\begin{proof}
Choose a smooth projective small model of \(X_0\). The transform of
\(M\) is big and movable. Run the projective minimal model program
for a small positive multiple of an effective rational representative
of this class. At each step the directing class remains movable.
A divisorial contraction would cover its exceptional prime by curves
of negative degree and would force that prime into the fixed part.
Bigness excludes a Mori fiber space. Thus the steps are flops of
smooth threefolds, and the basepoint-free theorem on the final model gives
\[
  g:T\longrightarrow Z,\qquad M_T=g^*H,
\]
with \(H\) ample on the normal projective variety \(Z\)
\cite{BCHM,KMbook,Kflops}. Since \(M_T\) is big, \(g\) is birational.
A volume form on \(T\) gives the zero canonical Weil divisor on \(Z\).
Thus \(K_Z\) is Cartier, \(g\) is crepant, and \(Z\) is canonical.

On a smooth projective resolution of the comparison write
\(p:W\to X_0\), \(q:W\to T\), and \(r=gq\).
For compatible divisors representing our rational classes,
\[
 p^*M=r^*H+G
\]
with \(G\) exceptional for both \(p\) and \(q\). For divisible \(m\),
strict transform identifies the complete systems on \(X_0\) and \(T\).
The pullbacks of their sections have divisors differing by \(mG\).
The system obtained from \(T\) is free, so \(mG\) is the fixed divisor
of the system obtained from \(X_0\). In particular \(G\geq0\).

If \(G\ne0\), the big negativity lemma
\cite[Proposition~9.4]{Graf}, applied to the projective birational
map \(r\), gives a component \(S\) on which \(-G|_S\) is big over
\(r(S)\). We can choose \(W\) initially so that its exceptional
divisors have simple normal crossings. Since \(S\) is \(q\)-exceptional,
\(\dim r(S)\leq1\). When \(r(S)\) is a point, \(S\) is a component of
the fiber by dimension. When \(r(S)\) is a curve, take a general
point away from the finitely many images of two-dimensional fibers.
A general fiber component \(V\) of \(S\to r(S)\) is then a component
of the full fiber of \(r\) and \(-G|_V\) is big. In the point case
put \(V=S\).
Write \(-G|_V\sim_{\QQ} A+B\) with \(A\) ample and \(B\) effective.
The full fiber is rationally chain connected
\cite[Corollary~1.6]{HM}. Hence a rational curve \(\Gamma\) in the
fiber passes through a general point of \(V\) outside \(B\) and the
other fiber components. It lies in \(V\) and has \(G\cdot\Gamma<0\).
Consequently \(p^*M\cdot\Gamma<0\), and \(p(\Gamma)\) is a rational
curve of negative \(M\)-degree. This contradicts the hypothesis.

It follows that \(p^*M=r^*H\). A sufficiently divisible multiple of
\(p^*M\) is generated, and its sections descend because
\(p_*\OO_W=\OO_{X_0}\). This proves semiampleness on \(X_0\).
The section rings agree, so the target is \(Z\), and the volume form
again gives crepancy.

\end{proof}

\begin{lemma}\label{lem:movable}
Let \(M\) be a rational line-bundle class on the family, with \(M_0\)
big and movable and \(M_t\) nef on a dense subset of the punctured
disc. Then \(M_0\) is semiample.
\end{lemma}
\begin{proof}
Since \(M_0\) is big and movable, some multiple has no fixed divisor. Its base locus has dimension at most one and hence only finitely many irreducible curve components. A rational curve of negative \(M_0\)-degree lies in this base locus, since a divisor in the system missing a general point of the curve would give it nonnegative degree. Lemma~\ref{lem:movement} produces infinitely many such curves of the same degree, a contradiction. Now apply Lemma~\ref{lem:semiample}.
\end{proof}

\subsection{The fixed part of a punctured polarization}

Choose \(L\in\Pic(\XX)\) with \(L_t\) ample on the punctured disc after
shrinking. The stable fixed primes of its big central restriction are
the components of the negative part of its divisorial Zariski
decomposition. They can be detected on any smooth projective
small model, since strict transform identifies prime divisors and
all complete section spaces.

\begin{lemma}\label{lem:fixed}
Let \(Y\) be a smooth projective small model of \(X_0\), and let
\(L_Y\) be the strict transform of \(L_0\). Each prime divisor
\(F\) in the stable fixed part of \(L_Y\) is an SQM-conic fibration
surface over a curve of positive genus in the sense of
\cite[Definitions~1.2 and~3.9]{Lutz}. For each such \(F\) there is
a smooth projective small model \(Y_F\) and a primitive divisorial
contraction of \(Y_F\) whose exceptional divisor is the strict transform
of \(F\) and whose restriction to it is its conic fibration over a smooth
curve of positive genus.
\end{lemma}

\begin{proof}
Write \(D=L_Y\). Proposition~\ref{prop:comparison} gives a smooth
projective comparison family \(\mathcal Y\to\DD\), with \(Y_0=Y\),
and small modifications \(\phi_t:X_t\dashrightarrow Y_t\).
Under the marking in that proposition, the class of \(D\) on
\(Y_t\) is the strict transform of \(L_t\). Fix a sufficiently
small \(t_*\ne0\). Strict transform identifies the real divisor class
spaces and carries movable classes to movable classes. It therefore
maps the open ample cone of \(X_{t_*}\) into the closed movable cone
of \(Y_{t_*}\). Since \(L_{t_*}\) is ample, its transform satisfies
\[
 [D]\in (\phi_{t_*})_*\operatorname{Amp}(X_{t_*})
 \subseteq \operatorname{int}\overline{\operatorname{Mov}}(Y_{t_*}).
\]

Let \(S\) be a small contractible representative of the Kuranishi
base of \(Y\), and mark the cohomology of its fibers by parallel
transport. The base is smooth by the Bogomolov--Tian--Todorov
theorem \cite{Tian,Todorov}. Since
\(H^1(Y,\OO_Y)=H^2(Y,\OO_Y)=0\), an ample line bundle on \(Y\)
extends over the family and remains ample after shrinking \(S\).
Write \(b\in S\) for the image of \(t_*\) under the classifying
map, so that the fiber over \(b\) is \(Y_{t_*}\).
After shrinking the Kuranishi base, the relative canonical bundle is
trivial. Indeed, $h^{3,0}=1$ is constant on these projective fibers,
so a generator on $Y$ extends locally by base change and is nowhere
zero near $Y$ by properness. Contraction with the volume form then identifies the holomorphic tangent
bundle of each fiber with its bundle of holomorphic two-forms.
Thus \(h^1(T_{Y_s})=h^{2,1}(Y_s)\) is constant and
\(H^0(Y_s,T_{Y_s})=0\). The Kodaira--Spencer map is an isomorphism
at the origin and hence throughout a neighborhood of it. Consequently
the family over a small neighborhood of \(b\) is also a Kuranishi
family for \(Y_{t_*}\).

Choose a point \(u\) in this neighborhood that is very general for both
the original and the recentered family. This is possible since the union
of their two countable collections of proper analytic jumping loci is
meagre. By \cite[Proposition~3.1]{Lutz}, the closed movable cone of
\(Y_u\) represents both very general cones in the same marking. Applying
\cite[Lemma~3.3(6)]{Lutz} at \(b\) gives
\[
 \overline{\operatorname{Mov}}(Y_{t_*})
 \subseteq \overline{\operatorname{Mov}}(Y_u)
 =\overline{\operatorname{Mov}}(Y^{\rm gen}).
\]
It follows that \([D]\) is in the interior of the last cone. Since
that cone lies in the pseudoeffective cone, \([D]\) is also big on
the very general fiber.

Set \(D_0=D\). Suppose that
\([D_j]\notin\overline{\operatorname{Mov}}(Y)\) while \([D_j]\)
is big and in the closed movable cone of the very general fiber.
By \cite[Definition~3.9 and Proposition~3.10]{Lutz} there is an
SQM-conic fibration surface \(E_j\) over a curve of positive genus
such that \(D_j\cdot\ell_{E_j}<0\). Here \(\ell_{E_j}\) is the
marked integral class of a whole general conic fiber. Put
\[
 a_j=-D_j\cdot\ell_{E_j},\qquad D_{j+1}=D_j-a_jE_j .
\]
The class of \(D_{j+1}\) is the reflection of \([D_j]\) across the hyperplane
\(\ell_{E_j}^{\perp}\). By \cite[Corollary~4.14]{Lutz}, it is
again big and in the closed movable cone of the very general fiber.
By \cite[Lemma~3.3(3)]{Lutz}, every class big on that fiber is big
on \(Y\), so each \(D_j\) is big on \(Y\). All \(D_j\) are integral
and each \(a_j\) is a positive integer. For a fixed ample integral
divisor \(A\) on \(Y\), we have \(A^2\cdot E_j>0\) and
\[
 0<A^2\cdot D_{j+1}
   = A^2\cdot D_j-a_j(A^2\cdot E_j)
   < A^2\cdot D_j .
\]
The positive integers \(A^2\cdot D_j\) can decrease only finitely
many times. Thus for some \(r\ge0\) the class \([D_r]\) is big and
belongs to \(\overline{\operatorname{Mov}}(Y)\). Summing the
recursion in \(\Pic(Y)\) gives
\[
 L_Y=D_r+G,\qquad G=\sum_{j=0}^{r-1}a_jE_j\ge0 .
\]

By the decomposition of the effective movable cone into nef cones of
small modifications \cite[\S2, formula~(1)]{Lutz}, the transform of
\(D_r\) is nef and big on some projective small modification.
The basepoint-free theorem gives a globally generated multiple there.
Strict transform then gives a multiple \(mD_r\) whose complete
linear system on \(Y\) has no fixed prime divisor.
For any prime divisor \(F\) outside \(\operatorname{Supp}G\), choose
\(H\in|mD_r|\) that does not contain \(F\). Then
\(H+mG\in|mL_Y|\) also avoids \(F\). Thus the support of the stable
fixed part of \(L_Y\) is contained in \(\operatorname{Supp}G\).
Each of its primes is consequently one of the \(E_j\).
For each such prime, \cite[Proposition~3.11]{Lutz} gives a sequence
of flops followed by a primitive divisorial contraction of its strict
transform along the conic ruling. The flops preserve smoothness by
\cite[Corollary~4.11]{Kflops}, and the image is the smooth base
curve of that ruling.
\end{proof}

\subsection{The projective contraction}\label{subsec:projcontraction}

The section ring of \(L_0\) is finitely generated. Indeed, its strict
transform \(L_Y\) on a smooth projective small model is big. Choose an
effective rational representative \(D\sim_{\QQ}L_Y\). For small
\(\epsilon>0\) rational, \((Y,\epsilon D)\) is klt and
\(K_Y+\epsilon D\sim_{\QQ}\epsilon L_Y\). Finite generation
\cite{BCHM} therefore gives a Veronese subring generated in degree one.
The fixed divisor in its \(k\)-th graded piece is \(k\) times the
fixed divisor of its generating system. This gives a decomposition
\(L_0=P_0+N\) with \(N\geq0\) rational, \(P_0\) big and movable, and
\[
 H^0(X_0,mP_0)\xrightarrow{\ \,\cdot s_{mN}\ \,}
 H^0(X_0,mL_0)
\]
an isomorphism for sufficiently divisible \(m\).
The coefficients of \(N\) are the asymptotic fixed multiplicities.

Run the movable minimal model program for the transform of \(P_0\)
as in Lemma~\ref{lem:semiample}. On the resulting smooth projective
model \(T\) its transform is nef, big and semiample. Write
\[
 g:T\longrightarrow Z,\qquad P_T=g^*H,
\]
where \(Z\) is normal and projective and \(H\) is ample and rational.
The volume form gives crepancy of \(g\) and canonical singularities
of \(Z\).

\begin{lemma}\label{lem:all}
Every component of \(N_T\) is exceptional for \(g\). Every exceptional prime of \(g\) has the conic presentation of Lemma~\ref{lem:fixed}.
\end{lemma}

\begin{proof}
Let \(E\) occur in \(N_T\) with coefficient \(a>0\). For rational \(0\le s\le a\), multiplication by the section of \(m(N_T-sE)\) identifies \(H^0(T,m(P_T+sE))\) with a subspace of \(H^0(T,mL_T)\). For sufficiently divisible \(m\) the latter has fixed part \(mN_T\), so dividing by the section of \(msE\) identifies \(H^0(T,m(P_T+sE))\) with \(H^0(T,mP_T)\). The volume is constant in \(s\), and differentiability of the volume gives \(P_T^2E=0\) \cite{BFJ}. Since \(P_T=g^*H\), the projection formula implies that \(E\) is exceptional.

Let \(G_T\) be the sum of all \(g\)-exceptional primes. Its strict transform \(G_0\) on \(X_0\) extends as a line-bundle class by Lemma~\ref{lem:family}. Choose a small positive rational \(\epsilon\) for which \(L+\epsilon G\)
is ample at one sufficiently close punctured parameter, using openness
of the ample cone on that fiber. After clearing denominators,
Lemma~\ref{lem:amplelocus} makes this class ample outside a closed
discrete subset of the disc. On \(T\) the stable fixed part of this perturbed class is precisely \(N_T+\epsilon G_T\). Indeed, for every effective integral exceptional divisor \(F\), normality gives \(g_*\OO_T(F)=\OO_Z\), so the projection formula gives the claimed fixed part for all divisible multiples. Now apply Lemma~\ref{lem:fixed} to the perturbed class.
\end{proof}

\begin{lemma}\label{lem:one}
The fibers of \(g\) have dimension at most one. On each exceptional
prime, every component of a general conic of the intrinsic ruling has
degree zero against \(P_T\).
\end{lemma}

\begin{proof}
Fix an exceptional prime \(E\) and a small model on which its
transform \(E'\) carries a conic fibration \(E'\to B\), where
\(B\) has positive genus. By \cite[Proposition~4.6]{Kflops}, the two threefolds
are isomorphic away from finite unions of rational curves and points.
The part omitted on \(E'\) maps to a finite set \(S\subset B\).
In fact an omitted rational curve contained in \(E'\) maps to a
point of \(B\), and one not contained in \(E'\) meets it in finitely
many points.

Let \(R\subset E\times B\) be the closure of the graph of the ruling
on the common open set, with projections \(a:R\to E\) and \(b:R\to B\).
At an omitted point of \(E\), the possible values of \(b\) lie in
\(S\). To see this, use a resolution of the graph of \(E\dashrightarrow
E'\). The points over an omitted point map to the omitted set of \(E'\),
since on the common open set the inverse is holomorphic and takes
values in the corresponding open set of \(E\). Thus \(a\) is finite.
For any rational curve \(C\subset E\), the set \(b(a^{-1}(C))\)
is finite. If \(C\) meets the common open set its strict transform is
rational and has constant image in \(B\), and the remaining values lie
in \(S\). If \(C\) is contained in the omitted set, all its values lie
in \(S\).

It follows that only finitely many rational curves on \(E\) pass through
a prescribed point or meet a prescribed rational curve \(C\subset E\).
Indeed, a curve meeting \(C\) is dominated by an irreducible component
of its inverse image in \(R\), vertical over \(B\). This component
meets \(a^{-1}(C)\) because its projection is proper and surjective
onto the curve. It must therefore be one of the finitely many curve
components of the fibers over \(b(a^{-1}(C))\). The same proof applies
with a point in place of \(C\).

Suppose a fiber of \(g\) contains a surface. Each surface component is
an exceptional prime as above. A rational curve not contained in a
given surface component meets it in only finitely many points.
Together with the finitely many one-dimensional components of the
fiber, the preceding finiteness shows that only finitely many rational
curves in the fiber meet a fixed rational curve, and only finitely many
pass through a fixed point. Induction on chain length gives at most
countably many curves reachable from a fixed point. Their union cannot
contain a surface component, by the Baire category theorem. This
contradicts rational chain connectedness of the fibers of the resolution
\(g\) of the canonical threefold \(Z\) \cite[Corollary~1.6]{HM}.

Every exceptional prime therefore maps to a curve. Through a general
point of the prime passes a contracted rational curve in it, by rational
chain connectedness of the one-dimensional fiber. On the common open
set this curve is vertical for \(E'\to B\). Outside \(S\) whole conics lie
in that open set. Taking a general point on each local component of a
general conic therefore shows that all of its components are contracted.
\end{proof}

\subsection{The mixed intersection and the original central fiber}
\begin{theorem}\label{thm:main}
Let \(\pi:\XX\to\DD\) be a proper holomorphic submersion with connected three-dimensional fibers. Suppose that \(X_t\) is a simply connected projective Calabi--Yau threefold for every \(t\ne0\). After shrinking the disc, let \(L\) be a line bundle whose restrictions to the punctured fibers are ample. Then the central fiber is a simply connected Moishezon manifold with
trivial canonical bundle. Write the divisorial Zariski decomposition of
the big central class as \(L_0=P_0+N\), and use the marking to extend
\(P_0\) to a class \(P\). Then \(P_t\) is nef for all sufficiently small
\(t\ne0\), and \(P_0\) is semiample. The morphism \(f:X_0\to Z\) defined by a divisible multiple of \(P_0\)
is birational and crepant, the target \(Z\) is canonical and projective,
and all fibers of \(f\) have dimension at most one.
\end{theorem}
\begin{proof}
Lemma~\ref{lem:family} gives the assertions on topology and the
canonical bundle. For a prime component \(E_0\) of \(N\), take a smooth
projective resolution \(j:\widetilde E\to E_0\hookrightarrow X_0\)
which also resolves the rational map to the base of its intrinsic ruling.
General conics avoid the exceptional sets of the small comparisons,
as in Lemma~\ref{lem:one}. Hence the class \(F\) of a full fiber on
\(\widetilde E\) satisfies \(j^*P_0\cdot F=0\). The class \(F\) is
nonzero, nef and isotropic. Its orthogonal complement is negative
semidefinite by the surface Hodge index theorem. Thus
\(P_0^2E_0=(j^*P_0)^2\leq0\).
On \(T\), Lemma~\ref{lem:all} gives \(g_*N_T=0\) as a divisor cycle.
The projection formula now yields
\begin{equation}\label{eq:sign}
 (P_0^2N)_{X_0}\leq0=(P_T^2N_T)_T.
\end{equation}

Extend the comparison with \(T\) by Proposition~\ref{prop:comparison},
and denote the resulting projective family by \(\rho:\TT\to\DD\).
Write \(\mathcal P,\mathcal N\) for the marked transforms of the relative
classes \(P,L-P\) on \(\TT\). For \(m\) sufficiently large and divisible,
\(mP_T=g^*(mH)\) is generated and
\[
 H^i(T,mP_T)=H^i(Z,mH)=0\quad (i>0),
\]
since \(Rg_*\OO_T=\OO_Z\). Cohomology and base change therefore
make \(m\mathcal P\) relatively generated near zero.
On each nearby \(T_t\) it has positive top self-intersection, by
constancy of \(P_T^3>0\), and so is big. It follows by strict transform
that \(P_t\) is big and movable on \(X_t\).

Fix such a punctured \(t\). Run a \(P_t\)-minimal model program on
\(X_t\) with scaling of \(L_t\), using a small effective klt boundary
proportional to \(P_t\) and a general representative of a multiple of
the ample scaling class. All steps are elementary threefold flops,
as in Lemma~\ref{lem:semiample}. Denote the strict transforms on the
source of the \(j\)-th flop by \(P_j,L_j,N_j=L_j-P_j\), and choose an
effective generator \(\gamma_j\) of its contracted ray.
The scaling equality gives
\[
 (P_j+\lambda_j L_j)\cdot\gamma_j=0,\qquad
 P_j\cdot\gamma_j<0,\qquad \lambda_j>0 .
\]
In particular \(L_j\cdot\gamma_j>0\) and \(N_j\cdot\gamma_j>0\).
The final model and \(T_t\) map to the same normal projective variety
\(\operatorname{Proj}R(P_t)\). On both models the transformed class of
\(P_t\) is pulled back from its tautological ample rational class, and
the pushforwards of the respective transforms of \(N_t\) agree as
divisor cycles. The projection formula therefore identifies their
mixed intersections.

Polarizing Wilson's cubic transformation formula
\cite[Proposition~2.1 and (2.1.1)]{WilsonFlops} along the finite
sequence gives
\[
 (\mathcal P_t^2\mathcal N_t)_{T_t}
     -(P_t^2N_t)_{X_t}
 =-\sum_j c_j(P_j\cdot\gamma_j)^2(N_j\cdot\gamma_j),\qquad c_j>0,
\]
where \(N_t=L_t-P_t\). The two terms on the left are constant in
their respective smooth families. Their difference is
\((P_T^2N_T)_T-(P_0^2N)_{X_0}\geq0\) by \eqref{eq:sign}.
Every flop would contribute a strictly negative term, so no flop
occurs and \(P_t\) is nef.

Lemma~\ref{lem:movable} now gives semiampleness of \(P_0\).
Its morphism \(f:X_0\to Z\) has the same target as \(g\), by equality
of the section rings, and is crepant by Lemma~\ref{lem:semiample}.
Corresponding prime divisors have the same image under the two morphisms
since they agree on the common open set. A surface in a fiber of \(f\)
would give a prime on \(T\) mapping to a point, contrary to
Lemma~\ref{lem:one}.

\end{proof}

\begin{corollary}\label{prop:rankone-direct}
Let $X_0$ be as in Lemma~\ref{lem:family}. Take a projective terminal
small model $X_0\dashrightarrow T$ as in
Proposition~\ref{prop:comparison}(1). Suppose the strict transform
$H_0$ of a very ample line bundle $A$ on $T$ extends to a line bundle
ample on a dense set of punctured fibers. Then $X_0\dashrightarrow T$
is an isomorphism. In particular, $X_0$ is projective if $b_2(X_0)=1$.
\end{corollary}
\begin{proof}
Here only the existence assertion in Proposition~\ref{prop:comparison}
is needed. The transform of $|A|$ has no fixed divisor, so every
$H_0$-negative curve belongs to one of finitely many curve components
of its base locus. Lemma~\ref{lem:movement} therefore excludes negative
rational curves. On a smooth projective common resolution $p:W\to X_0$,
$q:W\to T$, one has $p^*H_0=q^*A+G$, where $G$ is the effective
fixed divisor and is $q$-exceptional. The negativity argument in the
proof of Lemma~\ref{lem:semiample} applies directly with $r=q$:
it only uses that the target is klt, $W$ is smooth projective, and
the fibers of $r$ are rationally chain connected. If $G$ were nonzero,
a $q$-contracted rational curve of negative $G$-degree would give an
$H_0$-negative rational curve on $X_0$. Thus $G=0$ and $|H_0|$ gives
a small morphism $h:X_0\to T$ with $H_0=h^*A$. Any nontrivial
fiber contains a rational curve, by projecting a rational chain from
the corresponding $q$-fiber. Its $H_0$-degree is zero. By
Lemma~\ref{lem:movement} there would be infinitely many distinct
such curves, all contained in the compact exceptional locus of $h$,
which has dimension at most one. Hence $h$ is an isomorphism.

If $b_2=1$, let $L$ be the punctured polarization of
Lemma~\ref{lem:family}. Both $H_0$ and $L_0$ are big. Their pullbacks
to $W$ pair positively with the square of any ample divisor on $W$,
so their classes lie on the same rational ray in $H^2(X_0,\QQ)$.
The Picard identifications of Lemma~\ref{lem:family} show that
positive powers of their marked extensions agree after killing
torsion. The extension of $H_0$ is therefore ample on every
sufficiently close punctured fiber.
\end{proof}

\section{The analytic minimal model program near a fiber}\label{sec:local}

\begin{lemma}\label{lem:germ}
Let \(p:V\to W\) be a proper bimeromorphic morphism of normal complex
spaces. Suppose that \(R^1p_*\OO_V=R^2p_*\OO_V=0\) near \(w\in W\)
and that \(\dim p^{-1}(w)\le1\). Let \(C_1,\ldots,C_r\) be the
irreducible curve components of the reduced fiber. Then
\[
 \varinjlim_{U\ni w}\Pic(p^{-1}(U))
   \xrightarrow[\phantom{xxx}]{\ (\deg_{C_i})_i\ } \mathbb Z^r
\]
is an isomorphism. In particular, \(p\) is projective over a
sufficiently small Stein neighborhood of \(w\).
\end{lemma}

\begin{proof}
Stein neighborhoods are cofinal. On each sufficiently small such
neighborhood \(U\), Leray and Cartan's theorem~B give
\(H^j(p^{-1}(U),\OO_V)=0\) for \(j=1,2\). The exponential sequence
therefore identifies the Picard group with \(H^2(p^{-1}(U),\mathbb Z)\).
Proper base change for sheaves on topological spaces identifies the
direct limit with \(H^2(p^{-1}(w),\mathbb Z)\)
\cite[Tag~09V6]{Stacks}. The normalization of a compact reduced curve
changes it topologically only by identifications of finitely many
points. Thus evaluation on the fundamental classes of its irreducible
components identifies its second cohomology with \(\mathbb Z^r\).
Isolated points contribute nothing. The resulting coordinates of a
Chern class are precisely the multidegrees. A line bundle representing
\((1,\ldots,1)\) on a neighborhood is ample on the fiber by the degree
criterion for proper curves and invariance of ampleness under reduction.
Openness of fiberwise ampleness gives the last assertion
\cite[Proposition~1.4]{NakZar}.
\end{proof}

We first extend the contraction \(f:X_0\to Z\) of
Theorem~\ref{thm:main}. Write \(P_0=f^*H\) with \(H\) ample on \(Z\),
and use the extension \(P\) on \(\XX\) from that theorem. Choose \(m\)
sufficiently divisible that \(mP\) is integral, \(mH\) is very ample
and \(H^i(Z,\OO_Z(mH))=0\) for \(i>0\). Since \(Z\) is canonical and
\(f\) is a resolution, \(Rf_*\OO_{X_0}=\OO_Z\). Projection to \(Z\)
gives the same vanishing for \(mP_0\). Grauert's theorem and
properness, applied to the evaluation map, show that \(mP\) is
relatively generated and that its direct image commutes with base
change after shrinking \(\DD\). It defines a morphism
\[
 h:\XX\longrightarrow\mathcal P=\PP_{\DD}(\pi_*\OO_{\XX}(mP))
\]
whose central restriction is \(i\circ f\), with \(i:Z\hookrightarrow
\mathcal P_0\) the embedding given by \(mH\).

For completeness, let \(\mathcal A=h_*\OO_{\XX}\). The direct-image
sequence for multiplication by the disc parameter \(t\), together with
\(R^1h_{0*}\OO_{X_0}=0\), gives
\(R^1h_*\OO_{\XX}=0\) near \(\mathcal P_0\) by Nakayama's lemma.
Shrink the disc using properness of \(\mathcal P\to\DD\). The same
sequence for \(t-a\) gives
\(\mathcal A/(t-a)\mathcal A=(h_a)_*\OO_{X_a}\) for every \(a\).
The unit map \(\OO_{\mathcal P}\to\mathcal A\) is surjective on the
central fiber, since \((h_0)_*\OO_{X_0}=i_*\OO_Z\), and hence is
surjective after further shrinking. Its quotient defines a closed
subspace \(\mathcal Z\subset\mathcal P\), and the Stein factorization
gives \(F:\XX\to\mathcal Z\) with \(F_0=f\). The sheaf \(\mathcal A\)
is torsion-free over the smooth base, so \(\mathcal Z\) is flat over
the disc. It is normal, as is the Stein target of each \(h_a\), since
the corresponding source is normal. Each \(F_a\) is bimeromorphic:
the basepoint-free bundle \(mP_a\) has positive top self-intersection
by its value on \(X_0\), hence \(h_a\) is generically finite, and its
Stein morphism has connected fibers. The same argument at a general
point makes \(F\) bimeromorphic. The proper fiber dimension theorem
and compactness give \(\dim F^{-1}(u)\le1\) near \(Z\). Finally
\(R^if_*\OO_{X_0}=0\) for \(i>0\), so the multiplication sequence for
\(t\) and Nakayama give \(R^iF_*\OO_{\XX}=0\) near \(Z\). By
Lemma~\ref{lem:germ}, \(F\) is projective on a small Stein neighborhood
of each point of \(Z\).

\begin{lemma}\label{lem:localflop}
Let \(F:V\to U\) be a projective bimeromorphic morphism over a disc,
where \(U\) is a normal Stein space and the fibers of \(F\) have
dimension at most one. Fix \(z\in U_0\). Near \(F^{-1}(z)\), assume
that \(V\) is smooth of relative dimension three over the disc, that
\(K_V\sim0\), and that \(F_0\) is bimeromorphic with normal canonical
target. Let \(D\) be a line bundle positive on every \(F\)-contracted
curve over the punctured disc. Consider a \(D\)-negative elementary
step of a klt MMP over \(U\) around \(\{z\}\), with adjoint divisor
a positive rational multiple of \(D\) up to relative
\(\QQ\)-linear equivalence. After shrinking \(U\) about \(z\), the step
is a flop, is an isomorphism over the punctured disc, and induces an
isomorphism of the two central models over \(U_0\). The new model
is smooth over the disc near its fiber over \(z\).
\end{lemma}

\begin{proof}
Write \(q:V\to Y\) for the extremal contraction. Since \(F_0\)
resolves a canonical space, the multiplication sequence for \(t\)
gives \(R^iF_*\OO_V=0\) near \(z\) for \(i>0\), as above.
Numerical equivalence over \(\{z\}\) tests against line bundles on
neighborhoods of its inverse image \cite[\S4]{FujinoMMP}. By
Lemma~\ref{lem:germ} the components of the fiber give independent
generators of its cone of curves. Thus the chosen ray contracts exactly
one reduced irreducible curve \(C\) over \(z\). Put \(y=q(C)\).
The negative of the adjoint divisor is \(q\)-ample along the fibers
over \(z\) and remains so near them. All curves contracted there have
negative \(D\)-degree and hence belong to \(V_0\). The exceptional
locus lies in the exceptional locus of the bimeromorphic central map
and has dimension at most two. Thus \(q\) is small in the fourfold,
and the klt flip \(q^+:V^+\to Y\) is a flop since \(K_V\sim0\).
It is an isomorphism over the locus where \(q\) is an isomorphism,
in particular over the punctured disc.

Write \(Y_0\) and \(V_0^+\) initially for the Weil divisors of \(t\).
Smallness identifies their prime components and multiplicities with
those of \(V_0\). The canonical form on \(V\) gives trivial canonical
divisors on \(Y\) and \(V^+\), also by smallness. Thus the pullbacks of
\(K_Y+Y_0\) to the two models are \(K_V+V_0\) and \(K_{V^+}+V_0^+\).
The smooth pair \((V,V_0)\) is plt. Discrepancy comparison shows that
the other two pairs are plt as well. Their boundaries are normal, and
their ambient spaces are klt, hence Cohen--Macaulay. The Cartier
fibers of \(t\) are therefore Cohen--Macaulay and generically reduced,
so they coincide with these boundaries scheme-theoretically.
Adjunction gives crepant maps
\[
       V_0\xrightarrow{q_0}Y_0\xleftarrow{q_0^+}V_0^+ .
\]
In particular \(Y_0\) is canonical, since \(V_0\) is smooth.

We claim that a unique \(q_0\)-exceptional prime \(E\) meets the
fiber over \(y\), as a germ along that fiber. Theorem~9.1 of
\cite{FujinoMMP} applied to the negative contraction implies that
\(C\) is rational. At its normalization \(u:\PP^1\to V\) the space
of maps has dimension at least \(\chi(\PP^1,u^*T_V)=4\), as in
Lemma~\ref{lem:movement}. For nearby maps the images under \(q\) lie
in a Stein neighborhood of \(y\), so are points. Their \(D\)-degrees
remain negative, forcing their images into \(V_0\). The reduced
\(q\)-fiber at \(y\) is \(C\), by the ray description and
connectedness. A nearby map with image in this fiber factors through
its normalization. Constancy of the degree of a \(q\)-ample line
bundle makes this factor a degree-one map, so these maps form a germ
of dimension three. There is consequently an analytic arc through
\(u\) whose image under evaluation in \(Y\) is nonconstant.
The resulting curves lie in the exceptional locus of \(q_0\).
Only finitely many irreducible components of this locus meet its
compact fiber after shrinking about \(y\). The irreducible family
over the arc belongs to one of them, and the variation of its images
forces that component to have dimension two. It gives an exceptional
prime \(E\) containing \(C\).

For a line bundle \(B\) on a neighborhood of \(q_0^{-1}(y)\), let
\(d_B=B\cdot C\), and fix a \(q_0\)-ample bundle \(J\). The bundle
\(d_J B-d_B J\) has degree zero on \(C\) and is therefore trivial
after shrinking about \(y\), by Lemma~\ref{lem:germ} for the
resolution \(q_0\). It follows that degrees of \(B\) on all contracted
curves in this neighborhood are proportional to those of \(J\).
In particular \(E\cdot C<0\), for otherwise the effective
exceptional divisor \(E\) would be relatively nef near \(y\), in
contradiction to the negativity lemma \cite[Lemma~1.3]{WangNeg}.
Every nearby contracted curve then has negative degree against \(E\)
and is contained in \(E\). This proves uniqueness. It also gives
analytic \(\QQ\)-factoriality of \((Y_0,y)\). Indeed, given a germ of
a Weil divisor, add a rational multiple of \(E\) to its strict
transform on \(V_0\) to obtain degree zero on \(C\). A multiple is
principal over a smaller neighborhood by the same multidegree
isomorphism. Pushing down its meromorphic equation makes the
corresponding multiple of the original divisor principal.

The morphism \(q_0^+\) has a positive-dimensional fiber at \(y\).
Otherwise \(q^+\) would be finite and hence an isomorphism near
\(y\). A Cartier multiple of the flipped ample divisor would then
descend to \(Y\), whose pullback under the small morphism \(q\) has
degree zero on \(C\), a contradiction. If \(q_0^+\) were small near
\(y\), the pushforward of a relatively ample Cartier divisor would
be \(\QQ\)-Cartier by the preceding paragraph. Its pullback would
agree with the original divisor by smallness, again contradicting
positivity on the fiber. Hence \(q_0^+\) extracts a prime divisor.
Crepancy gives it discrepancy zero over \(Y_0\). On the smooth
crepant resolution \(V_0\), every further exceptional divisor has
strictly positive discrepancy. Thus the only exceptional valuation of
discrepancy zero near \(y\) is the valuation of \(E\). Both central
models extract precisely this valuation and are isomorphic in
codimension one over the germ.

Let \(A^+\) be a \(q_0^+\)-ample Cartier divisor and \(A\) its strict
transform on \(V_0\). The two relative section algebras agree, since
on each normal source the divisor conditions for a meromorphic section
are determined in codimension one. If \(A\cdot C<0\), degree
proportionality makes \(A\) negative on every contracted curve near
\(y\). The fibers of \(E\to q_0(E)\) then force the local sections
of each positive multiple of \(A\) to vanish along \(E\). This
contradicts relative generation of a sufficiently large multiple of
\(A^+\) at the generic points of its corresponding prime. If
\(A\cdot C=0\), then \(A\) is trivial over the germ, so the relative
Proj of the common section algebra is \(Y_0\), contradicting the
positive-dimensional fiber of \(q_0^+\). Hence \(A\cdot C>0\).
By openness of fiberwise ampleness, \(A\) is \(q_0\)-ample near \(y\).
The relative Proj of the common algebra now identifies the two central
models over this neighborhood.

Every other point of \(Y\) over \(z\) lies in the isomorphism locus
of \(q\) and \(q^+\). The central identifications agree there, since
they agree on the dense common open set. Properness of \(Y\to U\)
gives a single neighborhood of \(z\) on which they glue. The
resulting central model is smooth. The local ring of \(V^+\) at a
point on it has the nonzero-divisor \(t\) and regular quotient by \(t\).
It is regular and \(t\) is part of a regular system of parameters,
so \(V^+\) is smooth over the disc. The central identification
preserves the fiber dimension bound, and the punctured identification
preserves the positivity of the transform of \(D\). Thus these
conclusions apply at the following elementary step after shrinking
about \(z\).
\end{proof}

\section{Projectivity and K\"ahler specialization in dimension three}
\subsection{A polarization for the local models}\label{sec:closure}

Let \(f:X_0\to Z\) and \(g:T\to Z\) be the two birational contractions
constructed in Theorem~\ref{thm:main} and
Section~\ref{subsec:projcontraction}. Their section rings agree, so
\(X_0\dashrightarrow T\) is a bimeromorphic map over \(Z\).
We claim that it is an isomorphism away from the inverse image of a finite
subset \(\Sigma\subset Z\). Choose closed sets \(A\subset X_0\) and
\(A'\subset T\) of dimension at most one with identified complements.
Let \(p,q\) denote the projections of the reduced graph to \(X_0,T\).
For \(x\in A\), a positive-dimensional \(p\)-fiber projects under \(q\)
onto a positive-dimensional closed analytic subset of
\(A'\cap g^{-1}(f(x))\). Indeed, over \(T\setminus A'\) the inverse map
is defined and takes values in \(X_0\setminus A\). Hence \(f(x)\) is
the point image of a curve component of \(A'\).
The same argument applies to the other projection. Taking for \(\Sigma\)
the point images of all contracted curve components of \(A\) and \(A'\)
proves the claim, since finite bimeromorphic maps onto normal spaces are
isomorphisms.

Let \(B_i\) be the distinct curve images of the exceptional divisors
of \(f\). The positive-dimensional fibers outside their union lie over
finitely many points, being supported on curve components of the
exceptional locus. On each \(B_i\) remove the points of \(\Sigma\), the
singular points, and the intersections with the other \(B_j\).
After removing finitely many additional points, each resulting
\(B_i^\circ\) is a stratum of a Whitney stratification for the
projective map \(g\), adapted to its exceptional locus. Over this
stratum the same stratification applies to \(f\). The strata above
\(B_i^\circ\) map submersively to it, and the proper stratified family
of exceptional fibers is locally topologically trivial.

Fix \(b\in B_i^\circ\). Choose a local function \(h\) on \(Z\) whose
restriction to \(B_i^\circ\) is a coordinate vanishing at \(b\), and
put \(S_b=(h=0)\). The differential of \(h\circ f\) is nonzero along
\(f^{-1}(b)\) by submersivity on each stratum above \(B_i^\circ\).
Thus \(f^{-1}(S_b)\) is smooth near that fiber. After shrinking,
\(S_b\setminus\{b\}\) is smooth because there \(f\) is an isomorphism and
\(h\circ f\) is a submersion. The Cartier divisor \(S_b\) is
Cohen--Macaulay and regular in codimension one, hence normal.
Adjunction gives a crepant resolution of \(S_b\), which is therefore
the minimal resolution of a Du Val surface germ. Write
\(C_1,\ldots,C_r\) for the reduced exceptional curves and \(R_b\)
for the intersection lattice of the resolved ball.

Take an ample line bundle on \(T\) and let \(H_0\) be its strict transform.
Since the comparison is an isomorphism over \(b\), each \(H_0\cdot C_j\)
is positive. The nonzero coefficients of a root
\(\alpha=\sum_j n_j[C_j]\in R_b\) all have the same sign, whence
\[
 \langle c_1(H_0),\alpha\rangle=\sum_jn_j(H_0\cdot C_j)\ne0 .
\]
In particular, every root has nonzero image in global real homology.
The component classes are transported along \(B_i^\circ\) by the
proper stratified family. Locally they give homotopies of compact
cycles in \(X_0\), so their global homology classes are constant under
transport. After refining the stratum one may keep the dual graph constant:
the incidence points and the individual curves are strata after a
local finite labeling of the components. Each component has
self-intersection \(-2\), and off-diagonal intersections count common
points. Transport thus preserves the Cartan lattice and permutes its
finite root set. The finitely many strata \(B_i^\circ\) therefore yield only
finitely many proper root hyperplanes in \(H^2(X_0,\RR)\).

Under the marking of Lemma~\ref{lem:family}, the ample cone of a fixed
punctured fiber is open in \(H^2(X_0,\RR)\). Choose a rational point of
this cone outside the root hyperplanes. After clearing denominators,
extend it to a line bundle \(D\) on \(\XX\). By
Lemma~\ref{lem:amplelocus}, after shrinking about zero, \(D_t\) is
ample on every punctured fiber. We keep fixed the extension
\(F:\XX\to\mathcal Z\) of \(f\) constructed in
Section~\ref{sec:local}.

\begin{lemma}\label{lem:root}
Let \(U\) be a neighborhood of \(b\in B_i^\circ\) in \(\mathcal Z\)
and \(F':V'\to U\) a proper model smooth over the disc near the fiber
over \(b\). Assume that it is identified with the original model
over the puncture and that \(F'_0\) is identified with \(f\) over
\(U_0\), in each case over the indicated base. If \(D'\) is a line
bundle on \(V'\) corresponding to \(D\) under the punctured
identification, then \(D'\cdot C'\ne0\) for every reduced irreducible
curve \(C'\subset (F'_0)^{-1}(b)\).
\end{lemma}

\begin{proof}
Extend \(h\) to a holomorphic function near \(b\) on \(U\).
In a local embedding of \(U\), take a sufficiently small closed ball
centered at \(b\), with boundary transverse to the smooth part of the
central slice. Use \(h=0\) and this same ball on both models. Their
central inverse images are identified resolved balls. On each, the
differentials of \(h\) pulled back to the model and of the parameter
\(t\) are independent, since the central resolved slice is smooth
and the model is smooth over the disc. Over the boundary the central
maps are isomorphisms, so the transversality also holds upstairs.
These rank and boundary conditions persist near the compact central
inverse images. Properness, followed by shrinking the parameter disc,
thus gives two proper smooth families of oriented four-manifolds with
boundary, with submersive boundary maps.

Parallel transport to a common punctured slice and back gives an
integral isometry between the two central intersection lattices.
It respects the pairings with the restrictions of the line bundles:
on the punctured slice the bundles are identified, and on each
family the Chern class restricts from the total space. The class of
\(C'\) has square \(-2\), so its transport to the first central
slice is a root \(\alpha\in R_b\). Hence
\(D'\cdot C'=\langle c_1(D_0),\alpha\rangle\ne0\).
\end{proof}

\begin{theorem}\label{thm:projective-1}
Let \(\pi:\XX\to\DD\) be a proper holomorphic submersion with connected
three-dimensional fibers. If \(X_t\) is a simply connected projective
Calabi--Yau threefold for every \(t\ne0\), then so is \(X_0\).
\end{theorem}

\begin{proof}
Fix \(z\in Z\). By Lemma~\ref{lem:germ}, \(F\) is projective over a
Stein neighborhood \(U\) of \(z\). Choose an \(F\)-ample divisor \(A\).
After shrinking \(U\), the line bundles of \(D\) and \(D-A\) have
effective representatives. Their direct images have rank one by
bimeromorphicity, and Cartan's theorem~A gives a section nonzero at the
generic point. Thus \(D\) is relatively big. For an effective
representative \(D_{\mathrm{eff}}\), take a sufficiently small rational
\(\varepsilon>0\) so that \((V,\varepsilon D_{\mathrm{eff}})\) is klt
near \(F^{-1}(z)\), where \(V=F^{-1}(U)\). Choose a general member
of a sufficiently divisible relatively very ample linear system as
scaling divisor, divided by its degree. After reducing \(\varepsilon\)
if necessary, the scaled pair is klt and its adjoint divisor is ample
over \(z\). The boundary $\varepsilon D_{\mathrm{eff}}$ is relatively big, as is the
unscaled adjoint divisor since $K_V\sim0$. The smooth space $V$ is
$\QQ$-factorial over the chosen point. The point \(\{z\}\) is Stein compact with noetherian
ring of germs. Thus \cite[Theorem~1.7]{FujinoMMP} gives a terminating
local program with scaling and final model \(F':V'\to U\) on which
the transform \(D'\) is nef on the curves over \(z\). Neighborhoods
of \(z\) may be shrunk at each step.

By Lemma~\ref{lem:localflop}, each step is small, induces an isomorphism
over the puncture and identifies the central models over \(U_0\).
These identifications preserve the hypotheses for the next step:
smoothness over the disc near the central fiber, trivial canonical
bundle, the upper bound one on fiber dimensions, and positivity of
the directing class on contracted curves over the puncture.
In particular \(D'\) is Cartier after shrinking, and
Lemma~\ref{lem:root} applies to the final model above every point of
\(U_0\cap B_i^\circ\).

Suppose that \(D'\cdot C=0\) for an irreducible reduced curve
\(C\subset(F')^{-1}(z)\). The fibers of \(f\) are rationally chain
connected. To see this, compose a smooth projective modification of \(X_0\) with
\(f\), apply \cite[Corollary~1.6]{HM} to this projective birational
morphism with canonical target, and project the chains down. Since
the fiber has dimension at most one, joining two distinct general
points of \(C\) shows that \(C\) is rational. Write
\(\nu:\PP^1\to C\subset V'\) for its normalization.

Let \(M\) be the reduced germ of \(\Hom(\PP^1,V')\) at \(\nu\).
The deformation estimate of Lemma~\ref{lem:movement} gives
\(\dim_\nu M\ge\chi(\PP^1,\nu^*T_{V'})=4\). For a nearby
map \(v\), its composite with \(F'\) is constant since \(U\) is
Stein. Denote its value by \(e(v)\). Evaluation at any fixed point of
\(\PP^1\) shows that \(e\) is holomorphic. For an \(F'\)-ample bundle
\(J\), the degrees of \(v^*J\) and \(v^*D'\) are locally constant,
respectively \(J\cdot C>0\) and zero. Fix a marked point with image
a smooth point of \(C\) lying on no other fiber component. If
\(e(v)=z\), its evaluation near this point forces the image to be
\(C\). The degree against \(J\) forces the induced map to the
normalization of \(C\) to have degree one. These maps are
reparametrizations of \(\nu\) and form a three-dimensional germ.
The inverse parametrization is obtained by evaluation at three
points with distinct smooth images. Some irreducible component of
\(M\) therefore has nonconstant \(e\). Restrict to an analytic arc
through \(\nu\) on which \(e\) is nonconstant.

The image curves along the arc have degree zero against \(D'\).
They lie in the central model because \(D_t\) is ample on the
punctured fibers. Their images in \(U_0\) trace a nonconstant arc in
the positive-dimensional fiber locus of \(f\), hence meet some
\(B_i^\circ\), since the complement of their union in that locus is finite. At such a
point the reduced image curve \(C_s\) satisfies
\(0=\deg(v_s^*D')=\deg(v_s:\PP^1\to C_s)(D'\cdot C_s)\),
contrary to Lemma~\ref{lem:root}. Thus \(D'\) is positive on
every curve over \(z\). By the degree criterion it is ample on the
fiber, including its possible nonreduced structure, and openness
of ampleness makes it relatively ample over a neighborhood of \(z\).

It follows that the final model over this neighborhood is the relative
\(\operatorname{Proj}\) of the algebra of \(D'\). Since each
total-space modification is small, the codimension-one description of
sections identifies this algebra, on every open subset, with
\[
 \mathcal R|_U,\qquad
 \mathcal R=\bigoplus_{k\ge0}F_*\OO_{\XX}(kD).
\]
The identifications use restriction to the common isomorphism locus
and respect multiplication. They therefore glue the local models
canonically. By compactness of \(Z\) and properness over the disc, a
finite collection of these neighborhoods covers \(\mathcal Z\) after
shrinking the disc. A common divisible \(m\) makes the tautological
sheaf \(\OO(m)\) invertible and relatively ample on the glued model.
The local identifications of its central fiber with \(X_0\) agree on
overlaps, since they agree over the dense isomorphism locus of \(f\).
Consequently \(\OO(m)\) restricts to a global \(f\)-ample line
bundle on \(X_0\). Tensoring with the pullback of a high power of
an ample line bundle on the projective target \(Z\) gives an ample
line bundle on \(X_0\) \cite[Lemma~3.3]{FujinoPositive}.
Lemma~\ref{lem:family} supplies the other assertions.
\end{proof}

\subsection{Finite covers and the two cases of positive irregularity}
\label{sec:kahlerthree}

All families in this subsection are proper holomorphic submersions over a
disc, with connected compact fibers. We first record how to retain the
marking of a finite cover and how to recognize an elliptic-curve times K3
factor from its holomorphic forms.

\begin{lemma}\label{lem:globalcover}
A connected finite \'etale cover of one fiber of \(\pi:\XX\to\DD\)
extends to a finite \'etale map \(p:\mathcal Y\to\XX\) for which
\(q=\pi\circ p\) has connected fibers.
\end{lemma}
\begin{proof}
Ehresmann's theorem and contractibility of the disc give a differentiable
trivialization. The inclusion of the chosen fiber induces an isomorphism
on fundamental groups, so the finite-index subgroup defining the cover
determines a topological cover of \(\XX\). Pulling back holomorphic charts
makes the latter map finite \'etale and recovers the prescribed complex
structure on the original cover. In a differentiable trivialization all
the restrictions correspond to the same subgroup, so all fibers of \(q\)
are connected.
\end{proof}

\begin{lemma}\label{lem:cosymplecticproduct}
Let \(Y\) be a compact complex threefold diffeomorphic to \(E\times S_0\),
where \(E\) is an elliptic curve and \(S_0\) is a K3 surface. Suppose
that \(Y\) carries a closed holomorphic one-form \(\eta\) and a holomorphic
two-form \(\sigma\) such that \(\eta\wedge\sigma\) is
nowhere zero. Then \(Y\) is biholomorphic to the product of an elliptic
curve and a K3 surface.
\end{lemma}
\begin{proof}
The period group \(\Lambda\) of \(\eta\) is generated by two elements,
since \(H_1(Y,\mathbb Z)\simeq\mathbb Z^2\). If the generators lay
on a real line, a constant multiple of \(\eta\) would have real periods.
Its imaginary part would be the differential of a global pluriharmonic
function. Such a function is constant on \(Y\), and a holomorphic
one-form with zero imaginary part vanishes. Hence \(\Lambda\) is a
lattice. Integration gives a proper holomorphic submersion
\[
 a:Y\longrightarrow B=\mathbb C/\Lambda,\qquad a^*dz=\eta .
\]
The homomorphism \(a_*:\pi_1(Y)\to\pi_1(B)\) is surjective by the
definition of \(\Lambda\), hence is an isomorphism between groups
isomorphic to \(\mathbb Z^2\). The homotopy sequence of the fibration
shows that its fibers are connected and, since \(\pi_2(B)=0\), simply
connected. The restriction of \(\sigma\) to a fiber \(S\) trivializes
\(K_S\), so \(S\) is a K3 surface.

There is a unique holomorphic vector field \(V\) with
\(\iota_V\sigma=0\) and \(\eta(V)=1\). Its complex flow is complete by
compactness and lifts translations of \(B\). It trivializes the pullback
of \(Y\) to \(\mathbb C\), so \(Y\) is the suspension of a homomorphism
\(\rho:\Lambda\to\operatorname{Aut}(S)\). Since
\(\pi_2(B)=\pi_3(B)=0\), inclusion identifies \(\pi_2(S)\) with
\(\pi_2(Y)\). Under this identification the monodromy action is the
action of the deck group of the universal cover on \(\pi_2(Y)\), up to
inversion. This action is trivial for the differentiable manifold
\(E\times S_0\). The Hurewicz theorem and universal coefficients now
give \(\rho(\lambda)^*=1\) on \(H^2(S,\mathbb Z)\) for every \(\lambda\).
By the uniqueness assertion of the K3 Torelli theorem
\cite{LooijengaPeters}, \(\rho\) is trivial. Thus \(Y\simeq B\times S\).
\end{proof}

\begin{lemma}\label{lem:marked-limits}
Let \(q:\mathcal Y\to\DD\) be a smooth family whose punctured fibers are
products of elliptic curves and K3 surfaces. Then \(Y_0\) is also the
product of an elliptic curve and a K3 surface.
\end{lemma}
\begin{proof}
Fix \(s\ne0\). On the punctured disc the dimensions of the spaces of
holomorphic one-forms and top forms are both one. By Grauert's theorem
the maps from the fibers of \(q_*\Omega^1_{\mathcal Y/\DD}\) and
\(q_*K_{\mathcal Y/\DD}\) at \(s\) to the corresponding section spaces
on \(Y_s\) are isomorphisms. Cartan's theorem~A on \(\DD\) thus gives
global relative forms \(\eta\) and \(\Omega\) restricting to nonzero
forms in the two degrees on \(Y_s\). The relative exterior derivative of
\(\eta\) vanishes by holomorphic continuation from the punctured disc.

The zero divisor of \(\Omega\) misses \(Y_s\). Each of its irreducible
components maps properly onto a closed analytic subset of \(\DD\) and
thus maps to a point. A vertical prime divisor is an entire fiber, since
each fiber is smooth and connected. After shrinking at zero and dividing
by a power of the parameter, replacing \(\Omega\)
by \(t^{-m}\Omega\) for the appropriate integer \(m\ge0\), we obtain
a trivialization of \(K_{\mathcal Y/\DD}\). Likewise, dividing \(\eta\)
by its order of divisibility by \(t\) along \(Y_0\) gives a relatively
closed one-form whose restriction \(\eta_0\) is not identically zero.
It is nonzero as a section on every sufficiently close fiber.

Choose an Ehresmann marking of \(H_1(Y_t,\mathbb Z)\simeq\mathbb Z^2\).
The periods \(\lambda_1(t),\lambda_2(t)\) of \(\eta_t\) on a marked
basis depend smoothly on \(t\). They are linearly independent over
\(\mathbb R\) even at zero, by the pluriharmonic-function argument in
Lemma~\ref{lem:cosymplecticproduct}. Let
\(\ell_t:\mathbb R^2\to\mathbb C\) send the standard basis to these
periods, and set \(B_t=\mathbb C/\ell_t(\mathbb Z^2)\). Integration
from a smooth choice of base point gives a holomorphic map
\(a_t:Y_t\to B_t\). It is surjective, since it is a nonconstant
holomorphic map from a compact connected complex manifold to a curve.
In the smooth identifications \(B_t\simeq T^2=\mathbb R^2/\mathbb Z^2\)
given by the induced quotient isomorphisms \(\bar\ell_t^{-1}\), the maps assemble into a proper differentiable map
\[
 A:\mathcal Y\longrightarrow T^2\times\DD,\qquad
 A(y)=\bigl(\bar\ell_t^{-1}(a_t(y)),t\bigr),\quad t=q(y).
\]
Indeed, local path integrals can be taken along paths varying smoothly
under an Ehresmann trivialization, and changes of path give precisely
the periods. Properness follows from that of \(q\). For \(t\ne0\), every
nonzero holomorphic one-form on the product comes from its elliptic
factor, and integration modulo its full period lattice identifies
\(a_t\) with the projection, up to an isomorphism of elliptic curves.

Let \(b\) be any regular value of \(a_0\), viewed in \(T^2\).
Then \((b,0)\) is a regular value of \(A\), since the derivative is already
surjective onto the torus directions along \(Y_0\), and \(q\) is a
submersion. The critical values of \(A\) form a closed set by properness.
On a small ball about \((b,0)\), Ehresmann's theorem makes \(A\) a
differentiable fiber bundle. It follows by comparison with a fiber over
\(t\ne0\) that \(a_0^{-1}(b)\) is connected and diffeomorphic to a K3
surface. Adjunction and the triviality of \(K_{Y_0}\) show that its
canonical bundle is trivial. Thus every regular fiber is a K3 surface.
There exist regular values by Sard's theorem, and their union of fibers
is a nonempty open subset of \(Y_0\).

Upper semicontinuity of \(h^0(Y_t,T_{Y_t})\) gives a nonzero holomorphic
vector field \(V\) on \(Y_0\), since translations act on every punctured
fiber. The holomorphic function \(\eta_0(V)\) is constant. If it vanished,
\(V\) would be tangent to the fibers of \(a_0\). Its restriction to every
regular fiber would be zero, since a K3 surface has no nonzero
holomorphic vector fields. It would then vanish on a nonempty open set,
a contradiction. Taking \(\sigma=\iota_V\Omega_0\), we get
\[
 \eta_0\wedge\sigma=\eta_0(V)\,\Omega_0 .
\]
This is nowhere zero. Lemma~\ref{lem:cosymplecticproduct} now proves the
assertion.
\end{proof}

\begin{theorem}\label{thm:kahlerthree}
Let \(\pi:\XX\to\DD\) be a smooth family of threefolds whose punctured
fibers are K\"ahler with vanishing real first Chern class. Then \(X_0\)
is K\"ahler and has vanishing real first Chern class.
\end{theorem}
\begin{proof}
The restrictions of \(c_1(T_{\XX/\DD})\) define a flat section, so the
vanishing of the real first Chern class holds throughout the disc.
Choose a Beauville--Bogomolov cover of one punctured fiber and extend
it to \(p:\mathcal Y\to\XX\) by Lemma~\ref{lem:globalcover}. All punctured
fibers of \(q=\pi\circ p\) are K\"ahler with vanishing real first Chern
class. Their Hodge numbers are constant by Hodge decomposition and
semicontinuity.

If the chosen cover is a torus, every fiber is a torus by Catanese's
theorem on deformations in the large \cite[Theorem~2.1]{Catanese}.
Suppose next that the cover is a product of an elliptic curve and a K3
surface. On each punctured fiber \(Y_t\) one has
\(h^{1,0}=h^{2,0}=h^{3,0}=1\). A nonzero canonical form is nowhere zero,
since its zero divisor has real Chern class zero and a nonzero
effective divisor has positive volume against a K\"ahler form. Choose
nonzero holomorphic forms \(\eta_t,\sigma_t\) in degrees one and two.
They are closed, and their conjugates generate \(H^1(Y_t,\OO)\) and
\(H^2(Y_t,\OO)\) by Hodge symmetry. Trivializing \(K_{Y_t}\) identifies
their cup product with the perfect pairing from Serre duality. Hence
\(\eta_t\wedge\sigma_t\) is a nonzero canonical form, thus nowhere zero.
All fibers are differentiably identified with the chosen cover, so
Lemma~\ref{lem:cosymplecticproduct} shows that every punctured fiber is
a product. Lemma~\ref{lem:marked-limits} gives the same conclusion at zero.

If the chosen cover is a simply connected strict Calabi--Yau threefold,
all fibers are simply connected. Hodge constancy on the punctured disc
gives \(h^{1,0}=h^{2,0}=0\) and \(h^{3,0}=1\) there. The same zero-divisor
argument trivializes each punctured canonical bundle. Since every real
degree-two cohomology class is of type \((1,1)\), the open K\"ahler cone
contains a rational class and Kodaira's embedding theorem gives
projectivity. Theorem~\ref{thm:projective-1} now applies to \(q\).
In all three cases \(Y_0\) is K\"ahler. The trace of a K\"ahler form
along \(p_0\) is a K\"ahler form on \(X_0\).
\end{proof}

\section{Propagation from one fiber}\label{sec:propagation}

\begin{proposition}\label{prop:genericK}
Let \(q:\mathcal Y\to\DD\) be a smooth family of threefolds. If one
fiber \(Y_*\) is a complex torus or a product of an elliptic curve and a
K3 surface, every fiber is of the same respective type. If \(Y_*\) is a
simply connected strict Calabi--Yau threefold, the projective locus
contains the complement of a closed discrete subset.
\end{proposition}
\begin{proof}
The first assertion for tori is Catanese's theorem
\cite[Theorem~2.1]{Catanese}. For the other cases, take a connected
neighborhood of the chosen parameter on which all fibers are K\"ahler
by deformation openness. Hodge numbers are constant there. Grauert
base change at the chosen parameter and Cartan's theorem~A on the disc
allow any prescribed holomorphic \(p\)-form on \(Y_*\) to be extended to
a global relative holomorphic \(p\)-form. Each such extension is
relatively closed by holomorphic continuation from the K\"ahler
neighborhood.

In the product case extend a one-form and a two-form pulled back from
the two factors, choosing them nonzero. Their wedge has no zeros on
\(Y_*\). By Remmert's theorem the image of its zero set is a proper
analytic subset \(\Sigma\) of the disc, hence closed and discrete.
Outside \(\Sigma\), Lemma~\ref{lem:cosymplecticproduct} identifies each
fiber as a product. For every \(s\in\Sigma\), choose a subdisc meeting
\(\Sigma\) only at \(s\). Lemma~\ref{lem:marked-limits} applied there
gives the same conclusion at \(s\).

In the strict Calabi--Yau case, \(h^{2,0}=0\) on the K\"ahler
neighborhood, so every fiber there is projective by Kodaira's embedding
theorem. The family has uncountably many projective fibers. By
\cite[Corollary~A.5]{Rao}, the projective locus is dense analytic
Zariski open, which proves the assertion.
\end{proof}

\begin{theorem}\label{thm:global-1}
Let \(\pi:\XX\to\DD\) be a smooth family of threefolds. If one fiber is
K\"ahler with vanishing real first Chern class, every fiber is K\"ahler
with vanishing real first Chern class.
\end{theorem}
\begin{proof}
The relative tangent bundle gives the assertion about Chern classes.
Extend a Beauville--Bogomolov cover of the given fiber to
\(q:\mathcal Y\to\DD\) as above. By Proposition~\ref{prop:genericK}, the
covered fibers are K\"ahler outside a closed discrete set \(\Sigma\).
Each \(s\in\Sigma\) has a disc neighborhood on which
Theorem~\ref{thm:kahlerthree} applies with \(s\) as origin. Consequently
every fiber of \(q\) is K\"ahler. Taking traces under the fiberwise
finite \'etale maps proves the theorem.
\end{proof}

\section{Non-K\"ahler specializations in higher dimensions}\label{sec:counterexamples}

The partial Mukai flops of Oguiso and Schr\"oer
\cite[Propositions~6.3--6.4 and Theorem~6.5]{OguisoSchroer},
combined with Huybrechts's deformation of an elementary transformation
\cite[Theorem~3.4]{HuybrechtsFlop}, give the following examples.

\begin{theorem}\label{prop:examples-even}
For every \(m\ge2\) there is a proper holomorphic submersion over a disc
whose fibers over the punctured disc are irreducible holomorphic symplectic
manifolds of dimension \(2m\) and whose central fiber is a smooth
non-K\"ahler Moishezon manifold carrying a holomorphic symplectic form.
For every odd \(n\ge5\) there is such a submersion with \(n\)-dimensional
fibers, with K\"ahler punctured fibers and a non-K\"ahler Moishezon
central fiber, all with trivial canonical bundle.
\end{theorem}

\begin{proof}
Let \(A=E\times F\) be a product of elliptic curves. The fiber over
zero of the summation map \(A^{[m+1]}\to A\) is the projective
generalized Kummer variety \(Y=K_m(A)\), an irreducible holomorphic
symplectic manifold of dimension \(2m\) \cite[\S7]{BeauvilleKummer}.
For \(u\in F[m+1]\), the subschemes in \(Y\) contained in
\(E\times\{u\}\) form the linear system
\(P_u=|(m+1)0_E|\simeq\PP^m\).
These \((m+1)^2\) projective spaces are disjoint and Lagrangian.
Lines \(\ell_u\subset P_u\) have the same real homology class in \(Y\).
Indeed, the linear systems over all \(u\in F\) form a product family
in \(A^{[m+1]}\), and restriction
\(H^2(A^{[m+1]},\RR)\to H^2(Y,\RR)\) is surjective
\cite[proof of Proposition~6.4]{OguisoSchroer}.

Flop \(P=P_0\) alone and denote the resulting smooth compact manifold
by \(Y^+\). In the common blowup \(p:W\to Y\), \(q:W\to Y^+\),
the exceptional divisor \(G\) is the incidence variety in
\(P\times P^\vee\), and \(\OO_W(G)\) has degree \(-1\) on a line of
either ruling. For \(a^+\in H^2(Y^+,\RR)\) let \(a\in H^2(Y,\RR)\)
be the corresponding class on the common open set. These cohomology
groups identify by excision and the Thom isomorphism since the
removed submanifolds have complex codimension \(m\ge2\).
The difference \(q^*a^+-p^*a\) is a multiple of \([G]\).
Evaluation on a line in each ruling gives
\[
 a^+\cdot\ell^+=-a\cdot\ell_0,
 \qquad a^+\cdot\ell_u=a\cdot\ell_0\quad(u\ne0),
\]
where \(\ell^+\) is a line in the dual projective space \(P^\vee\).
Thus the effective cycle \(\ell^++\ell_u\) pairs to zero with every
real degree-two class. A K\"ahler form would integrate positively over
this cycle. The manifold \(Y^+\) is Moishezon by bimeromorphic
invariance. Its holomorphic symplectic form extends from the common
open set by Hartogs's theorem. Its exterior derivative vanishes by
analytic continuation, and its top exterior power is nowhere zero
since a nonempty zero divisor cannot be supported on \(P^\vee\).

Huybrechts constructs a one-parameter deformation \(\mathcal Y\) of
\(Y\) in which contraction of the Kodaira--Spencer class with the
symplectic form gives a K\"ahler class on \(Y\). Its restriction to
\(P\) is nonzero. Thus the normal sequence is the nonsplit Euler
extension \cite[Lemma~3.6]{HuybrechtsFlop}
\[
 0\longrightarrow\Omega^1_P\longrightarrow
 N_{P/\mathcal Y}\simeq\OO_P(-1)^{\oplus(m+1)}
 \longrightarrow\OO_P\longrightarrow0.
\]
Blowing up \(P\) in the total space and contracting the exceptional
\(P\times P^\vee\) by the other projection gives a smooth proper
family with central fiber \(Y^+\) and with the same punctured family
\cite[Theorem~3.4]{HuybrechtsFlop}. After shrinking the disc, its
punctured fibers are irreducible holomorphic symplectic by deformation
openness. Finally, taking the product with a fixed elliptic curve
gives the odd-dimensional examples. A K\"ahler form on their central
fiber would restrict to a K\"ahler form on a slice isomorphic to \(Y^+\).
\end{proof}

\end{document}